\documentclass{article}%
\usepackage{amsmath}
\usepackage{amsfonts}
\usepackage{amssymb}
\usepackage{graphicx}%
\providecommand{\U}[1]{\protect \rule{.1in}{.1in}}
\newtheorem{theorem}{Theorem}

\newtheorem{corollary}[theorem]{Corollary}

\newtheorem{definition}[theorem]{Definition}
\newtheorem{example}[theorem]{Example}

\newtheorem{lemma}[theorem]{Lemma}

\newtheorem{proposition}[theorem]{Proposition}
\newtheorem{remark}[theorem]{Remark}

\newenvironment{proof}[1][Proof]{\noindent \textbf{#1.} }{\  $\Box$}
\begin{document}

\title{$G$-L\'{e}vy Processes under Sublinear Expectations}
\author{Mingshang HU and Shige PENG\thanks{The author thanks the partial support from The National Basic
Research Program of China (973 Program) grant No. 2007CB814900
(Financial Risk). }\\Institute of Mathematics\\Shandong
University\\250100, Jinan, China\\peng@sdu.edu.cn}

\maketitle

\begin{abstract}
We introduce $G$-L\'{e}vy processes which develop the theory of
processes with independent and stationary increments under the
framework of sublinear expectations. We then obtain the
L\'{e}vy-Khintchine formula and the existence for $G$-L\'{e}vy
processes. We also introduce $G$-Poisson processes.

\end{abstract}

\medskip

\noindent{\footnotesize{\bf Keywords:\hspace{2mm}Sublinear
expectation, $G$-normal distribution, $G$-Brownian motion,
$G$-expectation, L\'{e}vy process, $G$-L\'{e}vy process, $G$-Poisson
process,
L\'{e}vy-Khintchine formula, L\'{e}vy-It\^{o} decomposition}}% At least three keywords

%\noindent\footnotesize{\bf MSC(2000):\hspace{2mm} 60H10, 60H05, 60H30, 60E05, 60E07, 62C05, 62D05}%Five letters(Mathematics Subject Classification 2000)
 \vspace{2mm}
\baselineskip 15pt
\renewcommand{\baselinestretch}{1.10}
\parindent=16pt  \parskip=2mm
\rm\normalsize\rm

\section{Introduction}

Distribution and independence are two important notions in the
theory of probability and statistics. These two notions were
introduced in \cite{Peng2005, Peng2004} under the framework of
sublinear expectations. Recently a new central limit theorem (CLT)
under sublinear expectations has been obtained in \cite{Peng2008}
based on a new i.i.d. assumption. The corresponding limit
distribution of the CLT is a $G$-normal distribution. This new type
of sublinear distributions was firstly introduced in
\cite{Peng2006a} (see also [22-25]) for a new type of $G$-Brownian
motion and the related calculus of It\^{o}'s type.

$G$-Brownian motion has a very rich and interesting new structure
which non-trivially generalizes the classical one. Briefly speaking
a $G$-Brownian motion is a continuous process with independent and
stationary increments under a given sublinear expectation. A very
interesting new phenomenon of $G$-Brownian motion is that its
quadratic process is also a continuous process with independent and
stationary increments, and thus can be still regarded as a
$G$-Brownian motion. A natural problem is how to develop the theory
of L\'{e}vy processes, i.e., processes with independent and
stationary increments but not necessarily continuous, under
sublinear expectations. In particular, how to define Poisson
processes under sublinear expectations.

The purpose of this paper is to study the distribution property,
i.e., L\'{e}vy-Khintchine formula, of a L\'{e}vy process under
sublinear expectations. The corresponding L\'{e}vy-It\^{o}
decomposition will be discussed in our forthcoming work. We
introduce $G$-L\'{e}vy processes for simplicity and obtain that the
corresponding distributions satisfy a new type of nonlinear
parabolic integro-partial differential equations. Conversely, we can
directly construct $G$-L\'{e}vy processes from these type of
equations. A specific case is $G$-Poisson processes. By comparison
with classical methods, our methods are more simple and direct.
Books on L\'{e}vy processes, e.g., \cite{Bert, Levy, Sato}, are
recommended for understanding the present results.

This paper is organized as follows: in Section 2, we recall some
important notions and results of sublinear expectations and
$G$-Brownian motions. In Section 3 we introduce $G$-L\'{e}vy
processes. We discuss the characterization of $G$-L\'{e}vy processes
in Section 4. In Section 5 we obtain the L\'{e}vy-Khintchine formula
for $G$-L\'{e}vy processes. The existence of $G$-L\'{e}vy processes
is given in Section 6. For reader's convenience we present some
basic results of this new type of nonlinear parabolic
integro-partial differential equations in the Appendix.

\section{Basic settings}

We present some preliminaries in the theory of sublinear
expectations and the related $G$-Brownian motions. More details of
this section can be found in [19-25].

\subsection{Sublinear expectation}

Let $\Omega$ be a given set and let $\mathcal{H}$ be a linear space of real
functions defined on $\Omega$ such that if $X_{1},\ldots,X_{n}\in \mathcal{H}$,
then $\varphi(X_{1},\cdots,X_{n})\in \mathcal{H}$ for each $\varphi \in
C_{Lip}(\mathbb{R}^{n})$, where $C_{Lip}(\mathbb{R}^{n})$ denotes the space of
Lipschitz functions.

\begin{remark}
In particular, all constants belong to $\mathcal{H}$ and $|X|,$ $X^{+},$
$X^{-}\in \mathcal{H}$ if $X\in \mathcal{H}$.
\end{remark}

Here we use $C_{Lip}(\mathbb{R}^{n})$ in our framework only for some
convenience of techniques. In fact our essential requirement is that
$\mathcal{H}$ contains all constants and, moreover, $X\in \mathcal{H}$ implies
$|X|\in \mathcal{H}$. In general, $C_{Lip}(\mathbb{R}^{n})$ can be replaced by
other spaces for specific problem. We list other two spaces used in this paper.

\begin{itemize}
\item $C_{b.Lip}(\mathbb{R}^{n})$: the space of bounded and Lipschitz functions;

\item $C_{b}^{k}(\mathbb{R}^{n})$: the space of bounded and $k$-time
continuously differentiable functions with bounded derivatives of all orders
less than or equal to $k$.
\end{itemize}

\begin{definition}
A sublinear expectation $\mathbb{\hat{E}}$ on $\mathcal{H}$ is a functional
$\mathbb{\hat{E}}:\mathcal{H}\rightarrow \mathbb{R}$ satisfying the following
properties: for all $X,$ $Y\in \mathcal{H}$, we have

\begin{description}
\item[(a)] Monotonicity: $\mathbb{\hat{E}}[X]\geq \mathbb{\hat{E}}[Y]$ if
$X\geq Y$.

\item[(b)] Constant preserving: $\mathbb{\hat{E}}[c]=c$ for $c\in \mathbb{R}$.

\item[(c)] Sub-additivity: $\mathbb{\hat{E}}[X]-\mathbb{\hat{E}}%
[Y]\leq \mathbb{\hat{E}}[X-Y]$.

\item[(d)] Positive homogeneity: $\mathbb{\hat{E}}[\lambda X]=\lambda
\mathbb{\hat{E}}[X]$ for $\lambda \geq0$.
\end{description}

The triple $(\Omega,\mathcal{H},\mathbb{\hat{E}})$ is called a sublinear
expectation space (compare with a probability space $(\Omega,\mathcal{F},P)$).
\end{definition}

\begin{remark}
If the inequality in (c) is equality, then $\mathbb{\hat{E}}$ is a linear
expectation on $\mathcal{H}$. We recall that the notion of the above sublinear
expectations was systematically introduced by Artzner, Delbaen, Eber and Heath
\cite{ADEH1, ADEH2}, in the case where $\Omega$ is a finite set, and by
Delbaen \cite{Delb} for the general situation with the notation of risk
measure: $\rho(X):=\mathbb{\hat{E}}[-X]$. See also Huber \cite{Huber} for even
earlier study of this notion $\mathbb{\hat{E}}$ (called the upper expectation
$\mathbf{E}^{\ast}$ in Ch. 10 of \cite{Huber}).
\end{remark}

Let $X=(X_{1},\ldots,X_{n})$, $X_{i}\in \mathcal{H}$, denoted by $X\in
\mathcal{H}^{n}$, be a given $n$-dimensional random vector on a sublinear
expectation space $(\Omega,\mathcal{H},\mathbb{\hat{E}})$. We define a
functional on $C_{Lip}(\mathbb{R}^{n})$ by
\[
\mathbb{\hat{F}}_{X}[\varphi]:=\mathbb{\hat{E}}[\varphi(X)]\text{ \ for all
}\varphi \in C_{Lip}(\mathbb{R}^{n}).
\]
The triple $(\mathbb{R}^{n},C_{Lip}(\mathbb{R}^{n}),\mathbb{\hat{F}}_{X}%
[\cdot])$ forms a sublinear expectation space. $\mathbb{\hat{F}}_{X}$ is
called the distribution of $X$.

\begin{remark}
If the distribution $\mathbb{\hat{F}}_{X}$ of $X\in \mathcal{H}$ is not a
linear expectation, then $X$ is said to have distributional uncertainty. The
distribution of $X$ has the following four typical parameters:%
\[
\bar{\mu}:=\hat{\mathbb{E}}[X],\  \  \underline{\mu}:=-\mathbb{\hat{E}%
}[-X],\  \  \  \  \  \  \  \  \bar{\sigma}^{2}:=\hat{\mathbb{E}}[X^{2}],\  \  \underline
{\sigma}^{2}:=-\hat{\mathbb{E}}[-X^{2}].
\]
The intervals $[\underline{\mu},\bar{\mu}]$ and $[\underline{\sigma}^{2}%
,\bar{\sigma}^{2}]$ characterize the mean-uncertainty and the
variance-uncertainty of $X$.
\end{remark}

The following simple properties are very useful in sublinear analysis.

\begin{proposition}
Let $X,$ $Y\in \mathcal{H}$ be such that $\hat{\mathbb{E}}[Y]=-\hat{\mathbb{E}%
}[-Y]$, i.e., $Y$ has no mean uncertainty. Then we have%
\[
\hat{\mathbb{E}}[X+Y]=\hat{\mathbb{E}}[X]+\hat{\mathbb{E}}[Y].
\]
In particular, if $\hat{\mathbb{E}}[Y]=\hat{\mathbb{E}}[-Y]=0$, then
$\hat{\mathbb{E}}[X+Y]=\hat{\mathbb{E}}[X]$.
\end{proposition}

\begin{proof}
It is simply because we have $\hat{\mathbb{E}}[X+Y]\leq \hat{\mathbb{E}%
}[X]+\hat{\mathbb{E}}[Y]$ and%
\[
\hat{\mathbb{E}}[X+Y]\geq \hat{\mathbb{E}}[X]-\hat{\mathbb{E}}[-Y]=\hat
{\mathbb{E}}[X]+\hat{\mathbb{E}}[Y].
\]

\end{proof}

Noting that $\hat{\mathbb{E}}[c]=-\hat{\mathbb{E}}[-c]=c$ for all
$c\in \mathbb{R}$, we immediately have the following corollary.

\begin{corollary}
For each $X\in \mathcal{H}$, we have $\hat{\mathbb{E}}[X+c]=\hat{\mathbb{E}%
}[X]+c$ for all $c\in \mathbb{R}$.
\end{corollary}

\begin{proposition}
For each $X,$ $Y\in \mathcal{H}$, we have%
\[
|\mathbb{\hat{E}}[X]-\mathbb{\hat{E}}[Y]|\leq \mathbb{\hat{E}}[X-Y]\vee
\mathbb{\hat{E}}[Y-X].
\]
In particular, $|\mathbb{\hat{E}}[X]-\mathbb{\hat{E}}[Y]|\leq \mathbb{\hat{E}%
}[|X-Y|]$.
\end{proposition}

\begin{proof}
By sub-additivity and monotonicity of $\mathbb{\hat{E}}[\cdot]$, it is easy to
prove the inequalities.
\end{proof}

We recall some important notions under sublinear expectations.

\begin{definition}
Let $X_{1}$ and $X_{2}$ be two $n$-dimensional random vectors defined
respectively on sublinear expectation spaces $(\Omega_{1},\mathcal{H}%
_{1},\mathbb{\hat{E}}_{1})$ and $(\Omega_{2},\mathcal{H}_{2},\mathbb{\hat{E}%
}_{2})$. They are called identically distributed, denoted by $X_{1}\overset
{d}{=}X_{2}$, if%
\[
\mathbb{\hat{E}}_{1}[\varphi(X_{1})]=\mathbb{\hat{E}}_{2}[\varphi
(X_{2})]\text{ \ for all }\varphi \in C_{Lip}(\mathbb{R}^{n}).
\]
It is clear that $X_{1}\overset{d}{=}X_{2}$ if and only if their distributions coincide.
\end{definition}

\begin{definition}
Let $(\Omega,\mathcal{H},\mathbb{\hat{E}})$ be a sublinear expectation space.
A random vector $Y=(Y_{1},\cdots,Y_{n})\in \mathcal{H}^{n}$ is said to be
independent from another random vector $X=(X_{1},\cdots,X_{m})\in
\mathcal{H}^{m}$ under $\mathbb{\hat{E}}[\cdot]$ if for each test function
$\varphi \in C_{Lip}(\mathbb{R}^{m}\times \mathbb{R}^{n})$ we have
\[
\mathbb{\hat{E}}[\varphi(X,Y)]=\mathbb{\hat{E}}[\mathbb{\hat{E}}%
[\varphi(x,Y)]_{x=X}].
\]
$\bar{X}=(\bar{X}_{1},\cdots,\bar{X}_{m})\in \mathcal{H}^{m}$ is said to be an
independent copy of $X$ if $\bar{X}\overset{d}{=}X$ and $\bar{X}$ is
independent from $X$.
\end{definition}

\begin{remark}
Under a sublinear expectation space $(\Omega,\mathcal{H},\mathbb{\hat{E}})$.
$Y$ is independent from $X$ means that the distributional uncertainty of $Y$
does not change after the realization of $X=x$. Or, in other words, the
\textquotedblleft conditional sublinear expectation\textquotedblright \ of $Y$
knowing $X$ is $\mathbb{\hat{E}}[\varphi(x,Y)]_{x=X}$. {In the case of linear
expectation, this notion of independence is just the classical one.}
\end{remark}

It is important to note that under sublinear expectations the condition
\textquotedblleft$Y$ is independent from $X$\textquotedblright \ does not imply
automatically that \textquotedblleft$X$ is independent from $Y$%
\textquotedblright. See the following example:

\begin{example}
We consider a case where $\mathbb{\hat{E}}$ is a sublinear expectation and
$X,Y\in \mathcal{H}$ are identically distributed with $\mathbb{\hat{E}%
}[X]=\mathbb{\hat{E}}[-X]=0$ and $\bar{\sigma}^{2}=\mathbb{\hat{E}}%
[X^{2}]>\underline{\sigma}^{2}=-\mathbb{\hat{E}}[-X^{2}]$. We also assume that
$\mathbb{\hat{E}}[|X|]=\mathbb{\hat{E}}[X^{+}+X^{-}]>0$, thus $\mathbb{\hat
{E}}[X^{+}]=\frac{1}{2}\mathbb{\hat{E}}[|X|+X]=$$\frac{1}{2}\mathbb{\hat{E}%
}[|X|]>0$. In the case where $Y$ is independent from $X$, we have%
\[
\mathbb{\hat{E}}[XY^{2}]=\mathbb{\hat{E}}[X^{+}\bar{\sigma}^{2}-X^{-}%
\underline{\sigma}^{2}]=(\bar{\sigma}^{2}-\underline{\sigma}^{2}%
)\mathbb{\hat{E}}[X^{+}]>0.
\]
But if $X$ is independent from $Y$ we have%
\[
\mathbb{\hat{E}}[XY^{2}]=0.
\]

\end{example}

\subsection{$G$-Brownian motion}

For a given positive integer $n$, we denote by $\langle x,y\rangle$ the scalar
product of $x$, $y\in \mathbb{R}^{n}$ and by $|x|=\langle x,x\rangle^{1/2}$ the
Euclidean norm of $x$. We also denote by $\mathbb{S}(d)$ the space of all
$d\times d$ symmetric matrices and by $\mathbb{R}^{n\times d}$ the space of
all $n\times d$ matrices. For $A$, $B\in \mathbb{S}(d)$, $A\geq B$ means that
$A-B$ is non-negative.

\begin{definition}
($G$-normal distribution with zero mean) A $d$-dimensional random vector
$X=(X_{1},\cdots,X_{d})$ on a sublinear expectation space $(\Omega
,\mathcal{H},\mathbb{\hat{E}})$ is said to be $G$-normally distributed if for
each $a\,$, $b\geq0$ we have
\[
aX+b\bar{X}\overset{d}{=}\sqrt{a^{2}+b^{2}}X,
\]
where $\bar{X}$ is an independent copy of $X$. Here the letter $G$ denotes the
function
\[
G(A):=\frac{1}{2}\mathbb{\hat{E}}[\langle AX,X\rangle]\text{ \  \ for }%
A\in \mathbb{S}(d).
\]

\end{definition}

It is easy to prove that $\mathbb{\hat{E}}[X_{i}]=\mathbb{\hat{E}}[-X_{i}]=0$
for $i=1,\ldots,d$ and the function $G$ is a monotonic and sublinear function.

Let $(\Omega,\mathcal{H},\mathbb{\hat{E}})$ be a sublinear expectation space,
$(X_{t})_{t\geq0}$ is called a $d$-dimensional process if $X_{t}\in
\mathcal{H}^{d}$ for each $t\geq0$.

\begin{definition}
($G$-Brownian motion) Let $G:$ $\mathbb{S}(d)\rightarrow \mathbb{R}$ be a given
monotonic and sublinear function. A $d$-dimensional process $(B_{t})_{t\geq0}$
on a sublinear expectation space $(\Omega,\mathcal{H},\mathbb{\hat{E}})$ is
called a $G$-Brownian motion if the following properties are satisfied:

\begin{description}
\item[(i)] $B_{0}=0$;

\item[(ii)] For each $t,s\geq0$, $B_{t+s}-B_{t}$ is independent from
$(B_{t_{1}},B_{t_{2}},\ldots,B_{t_{n}})$ for each $n\in \mathbb{N}$ and $0\leq
t_{1}\leq \cdots \leq t_{n}\leq t$;

\item[(iii)] $B_{t+s}-B_{t}\overset{d}{=}\sqrt{s}X$ for $t,s\geq0$, where $X$
is $G$-normally distributed.
\end{description}
\end{definition}

\begin{remark}
If $\mathbb{\hat{E}}$ is a linear expectation in the above two definitions,
then the function $G$ is a linear function, $X$ is classically normal and
$(B_{t})_{t\geq0}$ is classical Brownian motion.
\end{remark}

The above two definitions can be non-trivially generalized to the following situations.

\begin{definition}
($G$-normal distribution with mean uncertainty) A pair of $d$-dimensional
random vectors $(X,\eta)$ on a sublinear expectation space $(\Omega
,\mathcal{H},\mathbb{\hat{E}})$ is called $G$-distributed if for each $a\,$,
$b\geq0$ we have
\[
(aX+b\bar{X},a^{2}\eta+b^{2}\bar{\eta})\overset{d}{=}(\sqrt{a^{2}+b^{2}%
}X,(a^{2}+b^{2})\eta),
\]
where $(\bar{X},\bar{\eta})$ is an independent copy of $(X,\eta)$. Here the
letter $G$ denotes the function
\[
G(p,A):=\mathbb{\hat{E}}[\frac{1}{2}\langle AX,X\rangle+\langle p,\eta
\rangle]\text{ \  \ for }(p,A)\in \mathbb{R}^{d}\times \mathbb{S}(d).
\]

\end{definition}

Obviously, $X$ is $\bar{G}$-normally distributed with
$\bar{G}(A)=G(0,A)$. The distribution of $\eta$ can be seen as the
pure uncertainty of mean (see [22-25]). It is easy to prove that $G$
is a sublinear function monotonic in $A\in \mathbb{S}(d)$.

\begin{definition}
(generalized $G$-Brownian motion) Let $G:$ $\mathbb{R}^{d}\times
\mathbb{S}(d)\rightarrow \mathbb{R}$ be a given sublinear function monotonic in
$A\in \mathbb{S}(d)$. A $d$-dimensional process $(B_{t})_{t\geq0}$ on a
sublinear expectation space $(\Omega,\mathcal{H},\mathbb{\hat{E}})$ is called
a generalized $G$-Brownian motion if the following properties are satisfied:

\begin{description}
\item[(i)] $B_{0}=0$;

\item[(ii)] For each $t,s\geq0$, $B_{t+s}-B_{t}$ is independent from
$(B_{t_{1}},B_{t_{2}},\ldots,B_{t_{n}})$ for each $n\in \mathbb{N}$ and $0\leq
t_{1}\leq \cdots \leq t_{n}\leq t$;

\item[(iii)] $B_{t+s}-B_{t}\overset{d}{=}\sqrt{s}X+s\eta$ for $t,s\geq0$,
where $(X,\eta)$ is $G$-distributed.
\end{description}
\end{definition}

The construction of $G$-Brownian motion was first given in \cite{Peng2006a,
Peng2006b} and $G$-distributed random vector was given in \cite{Peng2008}.

Moreover, we have the characterization of the generalized $G$-Brownian motion
(see \cite{Peng2007, Peng2009}).

\begin{theorem}
\label{th1}Let $(X_{t})_{t\geq0}$ be a $d$-dimensional process defined on a
sublinear expectation space $(\Omega,\mathcal{H},\mathbb{\hat{E}})$ such that

\begin{description}
\item[(i)] $X_{0}=0$;

\item[(ii)] For each $t,s\geq0$, $X_{t+s}-X_{t}$ and $X_{s}$ are identically
distributed and $X_{t+s}-X_{t}$ is independent from $(X_{t_{1}},X_{t_{2}%
},\ldots,X_{t_{n}})$ for each $n\in \mathbb{N}$ and $0\leq t_{1}\leq \cdots \leq
t_{n}\leq t$;

\item[(iii)] $\lim_{t\downarrow0}\mathbb{\hat{E}}[|X_{t}|^{3}]t^{-1}=0$.
\end{description}

Then $(X_{t})_{t\geq0}$ is a generalized $G$-Brownian motion, where
\[
G(p,A)=\lim_{t\downarrow0}\mathbb{\hat{E}}[\frac{1}{2}\langle AX_{t}%
,X_{t}\rangle+\langle p,X_{t}\rangle]t^{-1} \  \  \  \ for \  \ (p,A)\in
\mathbb{R}^{d}\times \mathbb{S}(d).
\]

\end{theorem}

\begin{remark}
In fact, paths of $(X_{t})_{t\geq0}$ in the above theorem are continuous due
to the condition (iii) (see \cite{DHP, H-P}). In the following sections, we
consider L\'{e}vy processes without the condition (iii) which contain jumps.
\end{remark}

\section{$G$-L\'{e}vy processes}

A process $\{X_{t}(\omega):\omega \in \Omega,t\geq0\}$ defined on a sublinear
expectation space $(\Omega,\mathcal{H},\mathbb{\hat{E}})$ is called
c\`{a}dl\`{a}g if for each $\omega \in \Omega$, $\lim_{\delta \downarrow
0}X_{t+\delta}(\omega)=X_{t}(\omega)$ and $X_{t-}(\omega):=\lim_{\delta
\downarrow0}X_{t-\delta}(\omega)$ exists for all $t\geq0$. We now give the
definition of L\'{e}vy processes under sublinear expectations.

\begin{definition}
A $d$-dimensional c\`{a}dl\`{a}g process $(X_{t})_{t\geq0}$ defined on a
sublinear expectation space $(\Omega,\mathcal{H},\mathbb{\hat{E}})$ is called
a L\'{e}vy process if the following properties are satisfied:

\begin{description}
\item[(i)] $X_{0}=0$;

\item[(ii)] Independent increments: for each $t,$ $s>0$, the increment
$X_{t+s}-X_{t}$ is independent from $(X_{t_{1}},X_{t_{2}},\ldots,X_{t_{n}})$,
for each $n\in \mathbb{N}$ and $0\leq t_{1}\leq \cdots \leq t_{n}\leq t$;

\item[(iii)] Stationary increments: the distribution of $X_{t+s}-X_{t}$ does
not depend on $t$.
\end{description}
\end{definition}

\begin{remark}
If $(X_{t})_{t\geq0}$ is a L\'{e}vy process, then the finite dimensional
distribution of $(X_{t})_{t\geq0}$ is uniquely determined by the distribution
of $X_{t}$ for each $t\geq0$.
\end{remark}

\begin{proposition}
\label{pr2}Let $(X_{t})_{t\geq0}$ be a $d$-dimensional L\'{e}vy process
defined on a sublinear expectation space $(\Omega,\mathcal{H},\mathbb{\hat{E}%
})$. Then for each $A\in \mathbb{R}^{n\times d}$, $(AX_{t})_{t\geq0}$ is an
$n$-dimensional L\'{e}vy process.
\end{proposition}

\begin{proof}
By the definition of distribution and independence, it is easy to prove the result.
\end{proof}

Let $(X_{t})_{t\geq0}$ be a $d$-dimensional L\'{e}vy process defined on a
sublinear expectation space $(\Omega,\mathcal{H},\mathbb{\hat{E}})$. In this
paper, we suppose that there exists a $2d$-dimensional L\'{e}vy process
$(X_{t}^{c},X_{t}^{d})_{t\geq0}$ defined on a sublinear expectation space
$(\tilde{\Omega},\mathcal{\tilde{H}},\mathbb{\tilde{E}})$ such that the
distributions of $(X_{t}^{c}+X_{t}^{d})_{t\geq0}$ and $(X_{t})_{t\geq0}$ are
same. In this paper, we only consider the distribution property of
$(X_{t})_{t\geq0}$. Hence, we can suppose $X_{t}=X_{t}^{c}+X_{t}^{d}$ on the
same sublinear expectation space $(\Omega,\mathcal{H},\mathbb{\hat{E}})$.

\begin{remark}
In classical linear expectation case, by the L\'{e}vy-It\^{o} decomposition,
the above assumption of $(X_{t})_{t\geq0}$ obviously holds, where $(X_{t}%
^{c})_{t\geq0}$ is the continuous part and $(X_{t}^{d})_{t\geq0}$ is
the jump part.
\end{remark}

Furthermore, we suppose $(X_{t}^{c}+X_{t}^{d})_{t\geq0}$ satisfying the
following assumption:%
\begin{equation}
\lim_{t\downarrow0}\mathbb{\hat{E}}[|X_{t}^{c}|^{3}]t^{-1}=0;\text{
\  \ }\mathbb{\hat{E}}[|X_{t}^{d}|]\leq Ct\text{ \ for }t\geq0, \label{Main}%
\end{equation}
where $C$ is a constant.

\begin{remark}
By the assumption on $(X_{t}^{c})_{t\geq0}$, we know that $(X_{t}^{c}%
)_{t\geq0}$ is a generalized $G$-Brownian motion. The assumption on the jump
part $(X_{t}^{d})_{t\geq0}$ implies that it is of finite variation. The more
complicated situation will be discussed in our forthcoming work.
\end{remark}

\begin{example}
Suppose $(X_{t}^{d})_{t\geq0}$ is a $1$-dimensional positive L\'{e}vy process,
i.e., jumps are positive. Note that%
\[
\mathbb{\hat{E}}[X_{t+s}^{d}]=\mathbb{\hat{E}}[X_{t}^{d}]+\mathbb{\hat{E}%
}[X_{t+s}^{d}-X_{t}^{d}]=\mathbb{\hat{E}}[X_{t}^{d}]+\mathbb{\hat{E}}%
[X_{s}^{d}]
\]
and $\mathbb{\hat{E}}[X_{t}^{d}]$ is increasing in $t$, then we obtain
$\mathbb{\hat{E}}[X_{t}^{d}]=\mathbb{\hat{E}}[X_{1}^{d}]t$. Obviously, it
satisfies (\ref{Main}).
\end{example}

\begin{definition}
A $d$-dimensional L\'{e}vy process $(X_{t})_{t\geq0}$ is called a $G$-L\'{e}vy
process if there exists a decomposition $X_{t}=X_{t}^{c}+X_{t}^{d}$ for each
$t\geq0$, where $(X_{t}^{c},X_{t}^{d})_{t\geq0}$ is a $2d$-dimensional
L\'{e}vy process satisfying (\ref{Main}).
\end{definition}

By Proposition \ref{pr2}, We immediately have

\begin{proposition}
Let $(X_{t})_{t\geq0}$ be a $d$-dimensional $G$-L\'{e}vy process defined on a
sublinear expectation space $(\Omega,\mathcal{H},\mathbb{\hat{E}})$. Then for
each $A\in \mathbb{R}^{n\times d}$, $(AX_{t})_{t\geq0}$ is an $n$-dimensional
$G$-L\'{e}vy process.
\end{proposition}

\section{Characterization of $G$-L\'{e}vy processes}

Let $(X_{t})_{t\geq0}$ be a $d$-dimensional $G$-L\'{e}vy process
with the decomposition $X_{t}=X_{t}^{c}+X_{t}^{d}$. In this section,
we will show that for each given $\varphi \in
C_{b.Lip}(\mathbb{R}^{d}),$ $u(t,x):=\mathbb{\hat
{E}}[\varphi(x+X_{t})]$ is a viscosity solution of the following
equation:
\begin{equation}
\partial_{t}u(t,x)-G_{X}[u(t,x+\cdot)-u(t,x)]=0,\text{ }u(0,x)=\varphi(x),
\label{Main-eq1}%
\end{equation}

where $G_{X}[f(\cdot)]$ is a nonlocal operator defined by%
\begin{equation}
G_{X}[f(\cdot)]:=\lim_{\delta \downarrow0}\mathbb{\hat{E}}[f(X_{\delta}%
)]\delta^{-1}\text{ \  \ for }f\in C_{b}^{3}(\mathbb{R}^{d})\text{ with
}f(0)=0. \label{Nonlocal-oper}%
\end{equation}

We first show that the definition of $G_{X}[f(\cdot)]$ is meaningful. For this
we need the following lemmas.

\begin{lemma}
For each $\delta \leq1$, we have
\[
\mathbb{\hat{E}}[|X_{\delta}^{c}|^{p}]\leq C_{p}\delta^{p/2}\text{ \ for each
}p>0,
\]
where $C_{p}$ is a constant only depending on $p$.
\end{lemma}

\begin{proof}
This is a direct consequence of Theorem \ref{th1}.
\end{proof}

\begin{lemma}
\label{le4} For each given $f\in C_{b}^{3}(\mathbb{R}^{d})$ with $f(0)=0$, we
have%
\[
\mathbb{\hat{E}}[f(X_{\delta})]=\mathbb{\hat{E}}[f(X_{\delta}^{d})+\langle
Df(0),X_{\delta}^{c}\rangle+\frac{1}{2}\langle D^{2}f(0)X_{\delta}%
^{c},X_{\delta}^{c}\rangle]+o(\delta).
\]

\end{lemma}

\begin{proof}
It is easy to check that%
\begin{align*}
I_{\delta}^{1}  &  :=f(X_{\delta})-f(X_{\delta}^{c})-f(X_{\delta}^{d})\\
&  =\int_{0}^{1}\langle X_{\delta}^{c},Df(X_{\delta}^{d}+\alpha X_{\delta}%
^{c})-Df(\alpha X_{\delta}^{c})\rangle d\alpha
\end{align*}
and%
\begin{align*}
I_{\delta}^{2}  &  :=f(X_{\delta}^{c})-\langle Df(0),X_{\delta}^{c}%
\rangle-\frac{1}{2}\langle D^{2}f(0)X_{\delta}^{c},X_{\delta}^{c}\rangle \\
&  =\int_{0}^{1}\int_{0}^{1}\langle(D^{2}f(\alpha \beta X_{\delta}^{c}%
)-D^{2}f(0))X_{\delta}^{c},X_{\delta}^{c}\rangle \alpha d\beta d\alpha.
\end{align*}
Note that $Df$ is bounded, then we get%
\begin{align*}
\mathbb{\hat{E}}[|I_{\delta}^{1}|]  &  \leq(\mathbb{\hat{E}}[|X_{\delta}%
^{c}|^{3}])^{\frac{1}{3}}(\mathbb{\hat{E}}[\int_{0}^{1}|Df(X_{\delta}%
^{d}+\alpha X_{\delta}^{c})-Df(\alpha X_{\delta}^{c})|^{\frac{3}{2}}%
d\alpha])^{\frac{2}{3}}\\
&  \leq C(\mathbb{\hat{E}}[|X_{\delta}^{c}|^{3}])^{\frac{1}{3}}(\mathbb{\hat
{E}}[\int_{0}^{1}|Df(X_{\delta}^{d}+\alpha X_{\delta}^{c})-Df(\alpha
X_{\delta}^{c})|d\alpha])^{\frac{2}{3}}\\
&  \leq C_{1}(\mathbb{\hat{E}}[|X_{\delta}^{c}|^{3}])^{\frac{1}{3}%
}(\mathbb{\hat{E}}[|X_{\delta}^{d}|])^{\frac{2}{3}}\\
&  \leq C_{2}\delta^{7/6}=o(\delta).
\end{align*}
It is easy to obtain $\mathbb{\hat{E}}[|I_{\delta}^{2}|]\leq C\mathbb{\hat{E}%
}[|X_{\delta}^{c}|^{3}]=o(\delta).$ Noting that%
\[
|\mathbb{\hat{E}}[f(X_{\delta})]-\mathbb{\hat{E}}[f(X_{\delta}^{d})+\langle
Df(0),X_{\delta}^{c}\rangle+\frac{1}{2}\langle D^{2}f(0)X_{\delta}%
^{c},X_{\delta}^{c}\rangle]|\leq \mathbb{\hat{E}}[|I_{\delta}^{1}+I_{\delta
}^{2}|],
\]
we conclude the result.
\end{proof}

\begin{lemma}
\label{le3}Let $(p,A)\in \mathbb{R}^{d}\times \mathbb{S}(d)$ and $f\in C_{b}%
^{2}(\mathbb{R}^{d})$ with $f(0)=0$ be given. Then $\lim_{\delta \downarrow
0}\mathbb{\hat{E}}[f(X_{\delta}^{d})+\langle p,X_{\delta}^{c}\rangle+\frac
{1}{2}\langle AX_{\delta}^{c},X_{\delta}^{c}\rangle]\delta^{-1}$ exists.
\end{lemma}

\begin{proof}
We define%
\[
g(t)=\mathbb{\hat{E}}[f(X_{t}^{d})+\langle p,X_{t}^{c}\rangle+\frac{1}%
{2}\langle AX_{t}^{c},X_{t}^{c}\rangle].
\]
Obviously, $g(0)=0$. For each $t,s\in \lbrack0,1]$,%
\[
|g(t+s)-g(t)|\leq Cs+\mathbb{\hat{E}}[Y]\vee \mathbb{\hat{E}}[-Y]\leq C_{1}s,
\]
where $Y=\langle p+AX_{t}^{c},X_{t+s}^{c}-X_{t}^{c}\rangle+\frac{1}{2}\langle
A(X_{t+s}^{c}-X_{t}^{c}),X_{t+s}^{c}-X_{t}^{c}\rangle$. Thus $g(\cdot)$ is
differentiable almost everywhere on $[0,1]$. For each fixed $t_{0}<1$ such
that $g^{\prime}(t_{0})$ exists, we have%
\[
\frac{g(\delta)}{\delta}=\frac{g(t_{0}+\delta)-g(t_{0})}{\delta}%
-\Lambda_{\delta},
\]
where%
\begin{align*}
\Lambda_{\delta}  &  =\delta^{-1}(g(t_{0}+\delta)-\mathbb{\hat{E}}[f(X_{t_{0}%
}^{d})+f(X_{t_{0}+\delta}^{d}-X_{t_{0}}^{d})+\langle p,X_{t_{0}+\delta}%
^{c}\rangle \\
&  \qquad \qquad \qquad \qquad \qquad+\frac{1}{2}\langle AX_{t_{0}}^{c},X_{t_{0}%
}^{c}\rangle+\frac{1}{2}\langle A(X_{t_{0}+\delta}^{c}-X_{t_{0}}^{c}%
),X_{t_{0}+\delta}^{c}-X_{t_{0}}^{c}\rangle]).
\end{align*}
Similar to the above estimate, it is not difficult to prove that
$|\Lambda_{\delta}|\leq C\sqrt{t_{0}}$, where $C$ is a constant independent of
$\delta$ and $t_{0}$. Thus
\[
|\limsup_{\delta \downarrow0}g(\delta)\delta^{-1}-\liminf_{\delta \downarrow
0}g(\delta)\delta^{-1}|\leq2C\sqrt{t_{0}}.
\]
Letting $t_{0}\downarrow0$, we get the result.
\end{proof}

By the above two lemmas, we know that the definition of $G_{X}[f(\cdot)]$ is
meaningful. It is easy to check that $G_{X}[f(\cdot)]$ satisfies the following
properties: for each $f,$ $g\in C_{b}^{3}(\mathbb{R}^{d})$ with $f(0)=0,$
$g(0)=0,$

\begin{description}
\item[1)] Monotonicity: $G_{X}[f(\cdot)]\geq G_{X}[g(\cdot)]$ if $f\geq g$.

\item[2)] Sub-additivity: $G_{X}[f(\cdot)+g(\cdot)]\leq G_{X}[f(\cdot
)]+G_{X}[g(\cdot)].$

\item[3)] Positive homogeneity: $G_{X}[\lambda f(\cdot)]=\lambda G_{X}%
[f(\cdot)]$ for all $\lambda \geq0.$
\end{description}

Now we give the definition of viscosity solution for (\ref{Main-eq1}).

\begin{definition}
\label{de1} A bounded upper semicontinuous (lower semicontinuous) function $u$
is called a viscosity subsolution (viscosity supersolution) of the equation
(\ref{Main-eq1}) if $u(0,x)\leq \varphi(x)$ ($\geq \varphi(x)$) and for each
$(t,x)\in(0,\infty)\times \mathbb{R}^{d}$ and for each $\psi \in C_{b}^{2,3}$
such that $\psi \geq u$ ($\leq u$) and $\psi(t,x)=u(t,x)$, we have%
\[
\partial_{t}\psi(t,x)-G_{X}[\psi(t,x+\cdot)-\psi(t,x)]\leq0\text{ \ (}%
\geq0\text{)}.
\]
A bounded continuous function $u$ is called a viscosity solution of the
equation (\ref{Main-eq1}) if it is both a viscosity subsolution and a
viscosity supersolution.
\end{definition}

We now give the characterization of $G$-L\'{e}vy processes.

\begin{theorem}
\label{th6}Let $(X_{t})_{t\geq0}$ be a $d$-dimensional $G$-L\'{e}vy
process. For each $\varphi \in C_{b.Lip}(\mathbb{R}^{d})$, define
$u(t,x)=\mathbb{\hat {E}}[\varphi(x+X_{t})]$. Then $u$ is a
viscosity solution of the equation (\ref{Main-eq1}).
\end{theorem}

\begin{proof}
We first show that $u$ is a continuous function. Obviously,
$|u(t,x)-u(t,y)|\leq C|x-y|$. Note that%
\[
u(t+s,x)=\mathbb{\hat{E}}[\varphi(x+X_{s}+X_{t+s}-X_{s})]=\mathbb{\hat{E}%
}[u(t,x+X_{s})],
\]
then for $s\leq1$, $|u(t+s,x)-u(t,x)|\leq C\mathbb{\hat{E}}[|X_{s}^{c}%
+X_{s}^{d}|]\leq C_{1}\sqrt{s}$. Thus $u$ is continuous. For each fixed
$(t,x)\in(0,\infty)\times \mathbb{R}^{d}$ and $\psi \in C_{b}^{2,3}$ such that
$\psi \geq u$ and $\psi(t,x)=u(t,x)$, we have%
\[
\psi(t,x)=u(t,x)=\mathbb{\hat{E}}[u(t-\delta,x+X_{\delta})]\leq \mathbb{\hat
{E}}[\psi(t-\delta,x+X_{\delta})].
\]
Therefore,%
\begin{align*}
0  &  \leq \mathbb{\hat{E}}[\psi(t-\delta,x+X_{\delta})-\psi(t,x)]\\
&  =-\partial_{t}\psi(t,x)\delta+\mathbb{\hat{E}}[\psi(t,x+X_{\delta}%
)-\psi(t,x)+I_{\delta}]\\
&  \leq-\partial_{t}\psi(t,x)\delta+\mathbb{\hat{E}}[\psi(t,x+X_{\delta}%
)-\psi(t,x)]+\mathbb{\hat{E}}[|I_{\delta}|],
\end{align*}
where $I_{\delta}=\delta \int_{0}^{1}[\partial_{t}\psi(t,x)-\partial_{t}%
\psi(t-\beta \delta,x+X_{\delta})]d\beta$. It is easy to show that%
\[
\mathbb{\hat{E}}[|I_{\delta}|]\leq C\mathbb{\hat{E}}[\delta(\delta+|X_{\delta
}|)]=o(\delta).
\]
By the definition of $G_{X}$, we get%
\[
\partial_{t}\psi(t,x)-G_{X}[\psi(t,x+\cdot)-\psi(t,x)]\leq0.
\]
Hence, $u$ is a viscosity subsolution of (\ref{Main-eq1}). Similarly, we can
prove that $u$ is a viscosity supersolution of (\ref{Main-eq1}). Thus $u$ is a
viscosity solution of (\ref{Main-eq1}).
\end{proof}

\begin{remark}
We do not know the uniqueness of viscosity solution for (\ref{Main-eq1}). For
this, we need the following representation of $G_{X}$.
\end{remark}

\section{L\'{e}vy-Khintchine representation of $G_{X}$}

In this section, we give a representation of the infinitesimal generator
$G_{X}$, which can be seen as the L\'{e}vy-Khintchine formula for $G$-L\'{e}vy
processes. We first give some lemmas.

\begin{lemma}
Let $(p,A)\in \mathbb{R}^{d}\times \mathbb{S}(d)$ and $f\in C_{b.Lip}%
(\mathbb{R}^{d})$ with $f(0)=0$ and $f(x)=o(|x|)$ be given. Then $\lim
_{\delta \downarrow0}\mathbb{\hat{E}}[f(X_{\delta}^{d})+\langle p,X_{\delta
}^{c}\rangle+\frac{1}{2}\langle AX_{\delta}^{c},X_{\delta}^{c}\rangle
]\delta^{-1}$ exists.
\end{lemma}

\begin{proof}
Since $f(x)=o(|x|)$, there exists a sequence $\{ \delta_{n}:n\geq1\}$ such
that $\delta_{n}\downarrow0$ and $|f(x)|\leq \frac{1}{n}|x|$ on $|x|\leq
\delta_{n}$. For each fixed $\delta_{n}$, we can choose $f_{n}^{\varepsilon
}\in C_{b}^{2}(\mathbb{R}^{d})$ with $f_{n}^{\varepsilon}(0)=0$ such that%
\[
|f(x)-f_{n}^{\varepsilon}(x)|\leq \frac{4L\varepsilon}{\delta_{n}}|x|+\frac
{1}{n}|x|,
\]
where $L$ is the Lipschitz constant of $f$. Thus%
\[
\frac{\mathbb{\hat{E}}[|f(X_{\delta}^{d})-f_{n}^{\varepsilon}(X_{\delta}%
^{d})|]}{\delta}\leq(\frac{4L\varepsilon}{\delta_{n}}+\frac{1}{n}%
)\frac{\mathbb{\hat{E}}[|X_{\delta}^{d}|]}{\delta}\leq C(\frac{4L\varepsilon
}{\delta_{n}}+\frac{1}{n}).
\]
By Lemma \ref{le3} and the above estimate, we conclude the result by letting
first $\varepsilon \downarrow0$ and then $n\rightarrow \infty$.
\end{proof}

We denote by
\[
\mathcal{L}_{0}=\{f\in C_{b.Lip}(\mathbb{R}^{d}):f(0)=0\text{ and
}f(x)=o(|x|)\}
\]
and%
\[
\mathcal{L}=\{(f,p,q,A):f\in \mathcal{L}_{0}\text{ and }(p,q,A)\in
\mathbb{R}^{d}\times \mathbb{R}^{d}\times \mathbb{S}(d)\}.
\]
It is clear that $\mathcal{L}_{0}$ and $\mathcal{L}$ are both linear spaces.
Now we define a functional $\mathbb{\hat{F}}[\cdot]$ on $\mathcal{L}$ by%
\[
\mathbb{\hat{F}}[(f,p,q,A)]:=\lim_{\delta \downarrow0}\mathbb{\hat{E}%
}[f(X_{\delta}^{d})+\frac{\langle p,X_{\delta}^{d}\rangle}{1+|X_{\delta}%
^{d}|^{2}}+\langle q,X_{\delta}^{c}\rangle+\frac{1}{2}\langle AX_{\delta}%
^{c},X_{\delta}^{c}\rangle]\delta^{-1}.
\]
Similar to the proof of the above lemma, we know that the definition of
$\mathbb{\hat{F}}[\cdot]$ is meaningful.

\begin{lemma}
The functional $\mathbb{\hat{F}}:\mathcal{L}\rightarrow \mathbb{R}$ satisfies
the following properties:

\begin{description}
\item[(1)] $\mathbb{\hat{F}}[(f_{1},p,q,A_{1})]\geq \mathbb{\hat{F}}%
[(f_{2},p,q,A_{2})]$ if $f_{1}\geq f_{2}$ and $A_{1}\geq A_{2}$.

\item[(2)] $\mathbb{\hat{F}}[(f_{1}+f_{2},p_{1}+p_{2},q_{1}+q_{2},A_{1}%
+A_{2})]\leq \mathbb{\hat{F}}[(f_{1},p_{1},q_{1},A_{1})]+\mathbb{\hat{F}%
}[(f_{2},p_{2},q_{2},A_{2})].$

\item[(3)] $\mathbb{\hat{F}}[\lambda(f,p,q,A)]=\lambda \mathbb{\hat{F}%
}[(f,p,q,A)]$ \ for all $\lambda \geq0.$

\item[(4)] If $f_{n}\in \mathcal{L}_{0}$ satisfies $f_{n}\downarrow0$, then
$\mathbb{\hat{F}}[(f_{n},0,0,0)]\downarrow0$.
\end{description}
\end{lemma}

\begin{proof}
It is easy to prove (1), (2) and (3). We now prove (4). For each fixed
$0<\eta_{1}<\eta_{2}<\infty$, it is easy to check
\[
f_{n}(x)\leq(\sup_{0<|y|\leq \eta_{1}}\frac{f_{1}(y)}{|y|})|x|+(\sup_{\eta
_{1}\leq|y|\leq \eta_{2}}f_{n}(y))\frac{|x|}{\eta_{1}}+(\sup_{|y|\geq \eta_{2}%
}f_{1}(y))\frac{|x|}{\eta_{2}}.
\]
Thus%
\[
\mathbb{\hat{F}}[(f_{n},0,0,0)]\leq C(\sup_{0<|y|\leq \eta_{1}}\frac{f_{1}%
(y)}{|y|}+\frac{\sup_{\eta_{1}\leq|y|\leq \eta_{2}}f_{n}(y)}{\eta_{1}}%
+\frac{\sup_{y\in \mathbb{R}^{d}}f_{1}(y)}{\eta_{2}}).
\]
Noting that $\sup_{\eta_{1}\leq|y|\leq \eta_{2}}f_{n}(y)\downarrow0$, we have
\[
\lim_{n\rightarrow \infty}\mathbb{\hat{F}}[(f_{n},0,0,0)]\leq C(\sup
_{0<|y|\leq \eta_{1}}\frac{f_{1}(y)}{|y|}+\frac{\sup_{y\in \mathbb{R}^{d}}%
f_{1}(y)}{\eta_{2}}).
\]
Letting first $\eta_{1}\rightarrow0$ and then $\eta_{2}\rightarrow \infty$, we
obtain (4).
\end{proof}

By (2) and (3) of the above lemma, we immediately obtain that there exists a
family of linear functionals $\{F_{u}:u\in \mathcal{U}_{0}\}$ defined on
$\mathcal{L}$ such that%
\[
\mathbb{\hat{F}}[(f,p,q,A)]=\sup_{u\in \mathcal{U}_{0}}F_{u}[(f,p,q,A)].
\]
The proof can be found in \cite{Peng2008}. Note that (1) and (4) of the above
lemma, then for each $F_{u}$, by Daniell-Stone theorem, there exist
$(p^{\prime},q^{\prime},Q)\in$ $\mathbb{R}^{d}\times \mathbb{R}^{d}%
\times \mathbb{R}^{d\times d}$ and a unique measure $v$ on $(\mathbb{R}%
^{d}\backslash \{0\},\mathcal{B}(\mathbb{R}^{d}\backslash \{0\}))$ such that
\[
F_{u}[(f,p,q,A)]=\int_{\mathbb{R}^{d}\backslash
\{0\}}f(z)v(dz)+\langle
p,p^{\prime}\rangle+\langle q,q^{\prime}\rangle+\frac{1}{2}\text{\textrm{tr}%
}[AQQ^{T}].
\]
Thus%
\begin{equation}
\mathbb{\hat{F}}[(f,p,q,A)]=\sup_{(v,p^{\prime},q^{\prime},Q)\in
\mathcal{U}_{0}}\{ \int_{\mathbb{R}^{d}\backslash
\{0\}}f(z)v(dz)+\langle
p,p^{\prime}\rangle+\langle q,q^{\prime}\rangle+\frac{1}{2}\text{\textrm{tr}%
}[AQQ^{T}]\}. \label{Rep1}%
\end{equation}
In particular,%
\begin{equation}
\lim_{\delta \downarrow0}\mathbb{\hat{E}}[f(X_{\delta}^{d})]\delta^{-1}%
=\sup_{v\in \mathcal{V}}\int_{\mathbb{R}^{d}\backslash
\{0\}}f(z)v(dz),
\label{Rep2}%
\end{equation}
where%
\begin{equation}
\mathcal{V}=\{v:\exists(p^{\prime},q^{\prime},Q)\text{ such that }%
(v,p^{\prime},q^{\prime},Q)\in \mathcal{U}_{0}\}. \label{Meas}
\end{equation}
Taking specific $f$, we can easily prove that

\begin{description}
\item[(1)] For each $\varepsilon>0$, $\sup_{v\in \mathcal{V}}v(\{z:|z|\geq
\varepsilon \})<\infty.$

\item[(2)] For each $\varepsilon>0$, the restriction of $\mathcal{V}$ on the
set $\{z:|z|\geq \varepsilon \}$ is tight.

\item[(3)] $\sup_{v\in \mathcal{V}}\int_{\mathbb{R}^{d}}|z|v(dz)<\infty$.
\end{description}

In fact, it is easy to deduce that (3) implies (1) and (2). Similarly, it is
also easy to show that all $(p^{\prime},q^{\prime},Q)$ in $\mathcal{U}_{0}$
are bounded. Now we give the representation of $G_{X}$. For each $f\in
C_{b}^{3}(\mathbb{R}^{d})$ with $f(0)=0$, by Lemma \ref{le4} and the above
analysis, we have%
\begin{align*}
G_{X}[f(\cdot)]  &  =\sup_{(v,p^{\prime},q^{\prime},Q)\in
\mathcal{U}_{0}}\{ \int_{\mathbb{R}^{d}\backslash
\{0\}}(f(z)-\frac{\langle Df(0),z\rangle
}{1+|x|^{2}})v(dz)+\langle Df(0),p^{\prime}+q^{\prime}\rangle \\
&  \qquad \qquad \qquad \qquad \qquad \qquad \qquad \qquad \qquad \qquad+\frac{1}%
{2}\text{\textrm{tr}}[D^{2}f(0)QQ^{T}]\}.
\end{align*}
Note that $\sup_{v\in
\mathcal{V}}\int_{\mathbb{R}^{d}}|z|v(dz)<\infty$ , then
we have the following L\'{e}vy-Khintchine representation of $G_{X}$:%
\begin{equation}
G_{X}[f(\cdot)]=\sup_{(v,q,Q)\in \mathcal{U}}\{ \int_{\mathbb{R}^{d}%
\backslash \{0\}}f(z)v(dz)+\langle Df(0),q\rangle+\frac{1}{2}\text{\textrm{tr}%
}[D^{2}f(0)QQ^{T}]\}. \label{Main-rep}%
\end{equation}
We summarize the above discussions as a theorem.

\begin{theorem}
\label{th5}Let $(X_{t})_{t\geq0}$ be a $d$-dimensional $G$-L\'{e}vy
process. Then $G_{X}[f(\cdot)]$ has the L\'{e}vy-Khintchine
representation (\ref{Main-rep}), where $(q,Q)\in
\mathbb{R}^{d}\times \mathbb{R}^{d\times d}$ and $v$ is a
measure on $(\mathbb{R}^{d}\backslash \{0\},\mathcal{B}(\mathbb{R}%
^{d}\backslash \{0\}))$ satisfying%
\begin{equation}
\sup_{(v,q,Q)\in \mathcal{U}}\{ \int_{\mathbb{R}^{d}}%
|z|v(dz)+|q|+\text{\textrm{tr}}[QQ^{T}]\}<\infty. \label{Main-as}%
\end{equation}

\end{theorem}

We then immediately have the following theorem.

\begin{theorem}
Let $(X_{t})_{t\geq0}$ be a $d$-dimensional $G$-L\'{e}vy process. For each
$\varphi \in C_{b.Lip}(\mathbb{R}^{d})$, define $u(t,x)=\mathbb{\hat{E}%
}[\varphi(x+X_{t})]$. Then $u$ is the unique viscosity solution of the
following integro-partial differential equation:%
\begin{align}
\partial_{t}u(t,x)-\sup_{(v,q,Q)\in \mathcal{U}}\{ \int_{\mathbb{R}%
^{d}\backslash \{0\}}(u(t,x+z)  &  -u(t,x))v(dz)+\langle
Du(t,x),q\rangle
\nonumber \\
&  +\frac{1}{2}\text{\textrm{tr}}[D^{2}u(t,x)QQ^{T}]\}=0, \label{Main-eq}%
\end{align}
where $\mathcal{U}$ represents $G_{X}$.
\end{theorem}

\begin{proof}
By Theorem \ref{th5} and \ref{th6}, $u$ is a viscosity solution of
(\ref{Main-eq}). For the uniqueness, see appendix.
\end{proof}

\begin{remark}
The definition of viscosity solution for (\ref{Main-eq}) is the same as
Definition \ref{de1}.
\end{remark}

\section{Existence of $G$-L\'{e}vy processes}

We denote by $\Omega=\mathbb{D}_{0}(\mathbb{R}^{+},\mathbb{R}^{d})$
the space
of all $\mathbb{R}^{d}$-valued c\`{a}dl\`{a}g functions $(\omega_{t}%
)_{t\in \mathbb{R}^{+}}$, with $\omega_{0}=0$, equipped with the
Skorokhod topology. The corresponding canonical process is
$B_{t}(\omega)=\omega_{t}$
for $\omega \in \Omega$, $t\geq0$. We define%
\[
\mathcal{F}_{t}:=\sigma \{B_{s}:s\leq t\} \text{ and }\mathcal{F}%
=\bigvee_{t\geq0}\mathcal{F}_{t}.
\]

Following \cite{Peng2005, Peng2006a, Peng2006b}, for each fixed
$T\in
\lbrack0,\infty)$, we set%
\[
L_{ip}(\mathcal{F}_{T}):=\{ \varphi(B_{t_{1}\wedge
T},\ldots,B_{t_{n}\wedge T}):n\in \mathbb{N},t_{1},\ldots t_{n}\in
\lbrack0,\infty),\varphi \in C_{b.Lip}(\mathbb{R}^{d\times n})\}.
\]
It is clear that $L_{ip}(\mathcal{F}_{t})\subset
L_{ip}(\mathcal{F}_{T})$ for
$t\leq T$. We also set%
\[
L_{ip}(\mathcal{F}):=\bigcup_{n=1}^{\infty}L_{ip}(\mathcal{F}_{n}).
\]

Let $\mathcal{U}$ be given and satisfy (\ref{Main-as}). We consider
the corresponding integro-partial differential equation
(\ref{Main-eq}). For each given initial condition $\varphi \in
C_{b.Lip}(\mathbb{R}^{d})$, the viscosity solution $u^{\varphi}$ for
(\ref{Main-eq}) exists (see appendix). Furthermore, we have the
following theorem.

\begin{theorem}
Let $u^{\varphi}$ denote the viscosity solution of (\ref{Main-eq})
with the initial condition $\varphi \in C_{b.Lip}(\mathbb{R}^{d})$.
Then we have

\begin{description}
\item[1)] $u^{\varphi}\geq u^{\psi}$ if $\varphi \geq \psi$.

\item[2)] $u^{\varphi+\psi}\leq u^{\varphi}+u^{\psi}$.

\item[3)] $u^{\varphi+c}=u^{\varphi}+c$ for $c\in \mathbb{R}$.

\item[4)] $u^{\lambda \varphi}=\lambda u^{\varphi}$ for all $\lambda \geq0$.

\item[5)] $u^{\varphi}(t+s,x)=u^{u^{\varphi}(t,x+\cdot)}(s,0)$.
\end{description}
\end{theorem}

\begin{proof}
It is easy to check 3)-5). 1) and 2) are proved in appendix.
\end{proof}

We now introduce a sublinear expectation $\mathbb{\hat{E}}$ on $L_{ip}%
(\mathcal{F})$ by the following two steps:

Step 1. For each $\xi \in L_{ip}(\mathcal{F})$ of the form
$\xi=\varphi (B_{t+s}-B_{t})$, $t,s\geq0$ and $\varphi \in
C_{b.Lip}(\mathbb{R}^{d})$, we define
$\mathbb{\hat{E}}[\xi]=u(s,0)$, where $u$ is a viscosity solution of
(\ref{Main-eq}) with the initial condition $u(0,x)=\varphi(x)$.

Step 2. For each $\xi \in L_{ip}(\mathcal{F})$, we can find a $\phi
\in
C_{b.Lip}(\mathbb{R}^{d\times m})$ such that $\xi=\phi(B_{t_{1}},B_{t_{2}%
}-B_{t_{1}},\ldots,B_{t_{m}}-B_{t_{m-1}})$,
$t_{1}<t_{2}<\cdots<t_{m}$. Then
we define $\mathbb{\hat{E}}[\xi]=\mathbb{\phi}_{m}$, where $\mathbb{\phi}%
_{m}\in \mathbb{R}$ is obtained via the following procedure:%
\begin{align*}
\phi_{1}(x_{1},\ldots,x_{m-1})  &
=\mathbb{\hat{E}}[\phi(x_{1},\ldots
,x_{m-1},B_{t_{m}}-B_{t_{m-1}})];\\
\phi_{2}(x_{1},\ldots,x_{m-2})  &  =\mathbb{\hat{E}}[\phi_{1}(x_{1}%
,\ldots,x_{m-2},B_{t_{m-1}}-B_{t_{m-2}})];\\
&  \vdots \\
\phi_{m-1}(x_{1})  &  =\mathbb{\hat{E}}[\phi_{m-2}(x_{1},B_{t_{2}}-B_{t_{1}%
})];\\
\phi_{m}  &  =\mathbb{\hat{E}}[\phi_{m-1}(B_{t_{1}})].
\end{align*}
The related conditional expectation of $\xi$ under
$\mathcal{F}_{t_{j}}$ is
defined by%
\begin{align*}
\mathbb{\hat{E}}[\xi|\mathcal{F}_{t_{j}}]  &  =\mathbb{\hat{E}}[\phi(B_{t_{1}%
},B_{t_{2}}-B_{t_{1}},\ldots,B_{t_{m}}-B_{t_{m-1}})|\mathcal{F}_{t_{j}}]\\
&  =\phi_{m-j}(B_{t_{1}},\ldots,B_{t_{j}}-B_{t_{j-1}}).
\end{align*}

By the above theorem, it is easy to prove that
$\mathbb{\hat{E}}[\cdot]$ consistently defines a sublinear
expectation on $L_{ip}(\mathcal{F})$. Then $L_{ip}(\mathcal{F})$ can
be extended to a Banach space under the norm
$||X||:=\mathbb{\hat{E}}[|X|]$. We denote this space by $L_{G}^{1}%
(\mathcal{F})$. Note that $|\mathbb{\hat{E}}[X]-\mathbb{\hat{E}}%
[Y]|\leq \mathbb{\hat{E}}[|X-Y|]$, then $\mathbb{\hat{E}}[\cdot]$
can be extended as a continuous mapping on $L_{G}^{1}(\mathcal{F})$
which is still a sublinear expectation. Similarly, it is easy to
check that the conditional expectation
$\mathbb{\hat{E}}[\cdot|\mathcal{F}_{t}]$ can be also extended as
a continuous mapping $L_{G}^{1}(\mathcal{F})\rightarrow L_{G}^{1}%
(\mathcal{F}_{t})$. We now prove that the canonical process
$(B_{t})_{t\geq0}$ is a $G$-L\'{e}vy process. For this, we need the
following lemma.

\begin{lemma}
Let $\{X_{n}\}_{n=1}^{\infty}$ and $\{Y_{n}\}_{n=1}^{\infty}$ be two
sequences of $d$-dimensional random vectors on a sublinear
expectation space $(\Omega,\mathcal{H},\mathbb{\hat{E}})$. We assume
that $Y_{n}$ is independent from $X_{n}$ for $n=1,2,\ldots$. If
there exist $X$, $Y\in \mathcal{H}^{d}$
such that $\mathbb{\hat{E}}[|X_{n}-X|]\rightarrow0$ and $\mathbb{\hat{E}%
}[|Y_{n}-Y|]\rightarrow0$, then $Y$ is independent from $X$.
\end{lemma}

\begin{proof}
For each fixed $\varphi \in$ $C_{Lip}(\mathbb{R}^{2d})$, we define%
\[
\bar{\varphi}_{n}(x)=\mathbb{\hat{E}}[\varphi(x,Y_{n})]\text{ \ and \ }%
\bar{\varphi}(x)=\mathbb{\hat{E}}[\varphi(x,Y)].
\]
It is clear that $|\bar{\varphi}(x)-\bar{\varphi}_{n}(\bar{x})|\leq
C(|x-\bar{x}|+\mathbb{\hat{E}}[|Y_{n}-Y|])$. Thus $|\mathbb{\hat{E}}%
[\bar{\varphi}_{n}(X_{n})]-\mathbb{\hat{E}}[\bar{\varphi}(X)]|\leq
C(\mathbb{\hat{E}}[|X_{n}-X|]+\mathbb{\hat{E}}[|Y_{n}-Y|])$. Note
that
$\mathbb{\hat{E}}[\bar{\varphi}_{n}(X_{n})]=\mathbb{\hat{E}}[\varphi
(X_{n},Y_{n})]$, then we obtain $\mathbb{\hat{E}}[\bar{\varphi}%
(X)]=\mathbb{\hat{E}}[\varphi(X,Y)]$, which implies that $Y$ is
independent from $X$.
\end{proof}

\begin{theorem}
The canonical process $(B_{t})_{t\geq0}$ is a $G$-L\'{e}vy process.
\end{theorem}

\begin{proof}
Consider $\mathbb{D}_{0}(\mathbb{R}^{+},\mathbb{R}^{2d})$ and the
canonical process $(\bar{B}_{t},\tilde{B}_{t})_{t\geq0}$. Similar to
above, we can construct a sublinear expectation, still denoted by
$\mathbb{\hat{E}}[\cdot]$,
on $L_{ip}(\mathcal{F})$ via the following integro-pde:%
\begin{align*}
\partial_{t}w(t,x,y)-\sup_{(v,q,Q)\in \mathcal{U}}\{ \int_{\mathbb{R}%
^{d}\backslash \{0\}}(w(t,x,y+z)  &  -w(t,x,y))v(dz)+\langle D_{x}%
w(t,x,y),q\rangle \\
&  +\frac{1}{2}\text{\textrm{tr}}[D_{x}^{2}w(t,x,y)QQ^{T}]\}=0.
\end{align*}
It is easy to check that the distribution of
$(\bar{B}_{t})_{t\geq0}$
satisfies the following equation:%
\[
\partial_{t}u(t,x)-\sup_{(q,Q)\in \mathcal{U}}\{ \langle Du(t,x),q\rangle
+\frac{1}{2}\text{\textrm{tr}}[D^{2}u(t,x)QQ^{T}]\}=0.
\]
Thus $(\bar{B}_{t})_{t\geq0}$ is the generalized $G$-Brownian
motion. The
distribution of $(\tilde{B}_{t})_{t\geq0}$ satisfies the following equation:%
\[
\partial_{t}u(t,y)-\sup_{v\in \mathcal{V}}\int_{\mathbb{R}^{d}\backslash
\{0\}}(u(t,y+z)-u(t,y))v(dz)=0.
\]
We now show that $\tilde{B}_{t}$ belongs to $L_{G}^{1}(\mathcal{F})$
and
satisfies (\ref{Main}). Consider the function $\phi_{0}(y)=d+\sum_{i=1}%
^{d}y_{i}(\arctan y_{i}+\pi/2)$. Define
$\phi_{0}^{N}(y)=\phi_{0}(y-N)$ for $N>0$ and
\[
\phi_{i+1}^{N}(y)=\sup_{v\in
\mathcal{V}}\int_{\mathbb{R}^{d}\backslash
\{0\}}(\phi_{0}^{N}(y+z)-\phi_{0}^{N}(y))v(dz)\text{ \ for
}i=1,2,\ldots.
\]
Then it is not difficult to check that $u^{N}(t,y):=\sum_{i=0}^{\infty}%
\frac{t^{i}}{i!}\phi_{i}^{N}(y)$ is the solution of the above
equation with the initial condition $\phi_{0}^{N}$. It is easy to
check that $\phi _{0}(y)\geq \sum_{i=1}^{d}y_{i}^{+}$ and
$u^{N}(t,0)\rightarrow0$ as $N\rightarrow \infty$. Therefore, we
conclude that $|\tilde{B}_{t}|\in
L_{G}^{1}(\mathcal{F})$. Note that $u(t,y):=|y|+t\sup_{v\in \mathcal{V}}%
\int_{\mathbb{R}^{d}}|z|v(dz)$ is a viscosity supersolution of the
above equation, then (\ref{Main}) holds. It is also easy to check
that the distribution of $(\bar{B}_{t}+\tilde{B}_{t})_{t\geq0}$
satisfies (\ref{Main-eq}). By the above lemma, $(B_{t})_{t\geq0}$ is
a $G$-L\'{e}vy process.
\end{proof}

\begin{example}
We consider the following $1$-dimensional equation:%
\[
\partial_{t}u(t,x)-G_{\lambda}(u(t,x+1)-u(t,x))=0,\text{ }u(0,x)=\varphi(x),
\]
where $G_{\lambda}(a)=a^{+}-\lambda a^{-}$, $\lambda \in
\lbrack0,1]$. This equation is a special case of the above equation
with $\mathcal{V}=\{ \delta_{l}:l\in \lbrack \lambda,1]\}$. Thus we
can construct the corresponding sublinear expectation
$\mathbb{\hat{E}}[\cdot]$. The canonical process $(B_{t})_{t\geq0}$
is called the $G$-Poisson process under this sublinear expectation
$\mathbb{\hat{E}}[\cdot]$. We also have

\begin{itemize}
\item If $\varphi$ is increasing, then $\mathbb{\hat{E}}[\varphi
(x+B_{t})]=\sum_{i=0}^{\infty}\frac{t^{i}}{i!}\varphi(x+i)e^{-t}$.

\item If $\varphi$ is decreasing, then $\mathbb{\hat{E}}[\varphi
(x+B_{t})]=\sum_{i=0}^{\infty}\frac{(\lambda t)^{i}}{i!}\varphi
(x+i)e^{-\lambda t}$.
\end{itemize}

In particular, $\mathbb{\hat{E}}[B_{t}]=t$ and $-\mathbb{\hat{E}}%
[-B_{t}]=\lambda t$. Thus it characterizes a Poisson process with
intensity uncertainty in $[\lambda,1]$.
\end{example}

\begin{remark}
We consider the following integro-pde:%
\begin{align}
\partial_{t}u(t,x)-\sup_{(v,q,Q)\in \mathcal{U}}\{& \int_{\mathbb{R}%
^{d}\backslash \{0\}}(u(t,x+z)  -u(t,x)-\frac{\langle
Du(t,x),z\rangle
}{1+|z|^{2}})v(dz)\nonumber \\
&+\langle
Du(t,x),q\rangle+\frac{1}{2}\text{\textrm{tr}}[D^{2}u(t,x)QQ^{T}]\}
  =0,\label{Main-geq}%
\end{align}
where%
\[
\sup_{(v,q,Q)\in \mathcal{U}}\{ \int_{\mathbb{R}^{d}}(|z|I_{[|z|\geq1]}%
+|z|^{2}I_{[|z|<1]})v(dz)+|q|+\text{\textrm{tr}}[QQ^{T}]\}<\infty
\]
and%
\[
\lim_{\kappa \downarrow0}\sup_{v\in \mathcal{V}}\int_{|z|\leq \kappa}%
|z|^{2}v(dz)=0.
\]
For each given initial condition $\varphi \in
C_{b.Lip}(\mathbb{R}^{d})$, the viscosity solution $u^{\varphi}$ for
(\ref{Main-geq}) exists (see appendix). Thus we can construct the
corresponding sublinear expectation. Obviously, the canonical
process $(B_{t})_{t\geq0}$ is a L\'{e}vy process and has the
decomposition $B_{t}=B_{t}^{c}+B_{t}^{d}$. But it does not satisfy
(\ref{Main}) if $\sup_{v\in
\mathcal{V}}\int_{|z|\leq1}|z|v(dz)=\infty$.
\end{remark}

\section{Appendix}

In the appendix, we mainly consider the domination of viscosity
solutions for (\ref{Main-eq}) and (\ref{Main-geq}). We refer to
\cite{A-T, B-I, Cr, CIL, J-K, Peng2008} and the references therein.
For simplicity, we consider
the following type of integro-pde:%
\begin{equation}
\partial _{t}u(t,x)-G(Du(t,x),D^{2}u(t,x),u(t,x+\cdot ))=0,\ u(0,x)=\varphi
\in C_{b.Lip}(\mathbb{R}^{d}),  \label{eq1}
\end{equation}%
where $G:\mathbb{R}^{d}\times \mathbb{S}(d)\times
C_{b}^{1,2}\rightarrow \mathbb{R}$, $(t,x)\in Q_{T}:=(0,T)\times
\mathbb{R}^{d}$. We suppose $G$ satisfies the following condition:

\begin{description}
\item[(A1)] If $p_{k}\rightarrow p$, $X_{k}\rightarrow X$ and for $%
(t_{k},x_{k})\rightarrow (t,x),$ $\phi _{k}(t_{k},x_{k}+\cdot
)\rightarrow \phi (t,x+\cdot )$ locally uniform on $\mathbb{R}^{d}$,
$\phi _{k}$ uniformly bounded and $D^{n}\phi _{k}\rightarrow
D^{n}\phi $, $n=1,2$,
locally uniform on $Q_{T}$, then%
\[
G(p_{k},X_{k},\phi _{k}(t_{k},x_{k}+\cdot ))\rightarrow G(p,X,\phi
(t,x+\cdot )).
\]

\item[(A2)] If $X\geq Y$ and $(\phi -\psi )(t,\cdot )$ has a global minimum
at $x$, then
\[
G(p,X,\phi (t,x+\cdot ))\geq G(p,Y,\psi (t,x+\cdot )).
\]

\item[(A3)] For each constant $c\in \mathbb{R}$, $G(p,X,\phi (t,x+\cdot
)+c)=G(p,X,\phi (t,x+\cdot ))$.
\end{description}

For the equation (\ref{Main-eq}), the corresponding $G$ is%
\[
G(p,X,u(t,x+\cdot ))=\! \! \! \! \sup_{(v,q,Q)\in \mathcal{U}}\! \! \{ \int_{%
\mathbb{R}^{d}\backslash \{0\}}\! \! \!
\!(u(t,x+z)-u(t,x))v(dz)+\langle p,q\rangle
+\frac{1}{2}\text{\textrm{tr}}[XQQ^{T}]\}.
\]%
Obviously, it satisfies all the above assumptions. The above
assumptions also hold for the equation (\ref{Main-geq}).

The definition of viscosity solution for (\ref{eq1}) is the same as
Definition \ref{de1}. We also suppose that for each given $\kappa \in (0,1)$%
, there exists $G^{\kappa }:\mathbb{R}^{d}\times \mathbb{S}(d)\times \mathrm{%
SC}_{b}(Q_{T})\times C^{1,2}(Q_{T})\rightarrow \mathbb{R}$, where $\mathrm{SC%
}_{b}(Q_{T})$ denotes the set of bounded upper or lower
semicontinuous
functions, satisfying the following assumptions: for $p\in \mathbb{R}^{d}$, $%
X$, $Y\in \mathbb{S}(d)$, $u$, $-v\in \mathrm{USC}_{b}(Q_{T})$,
$w\in \mathrm{SC}_{b}(Q_{T})$, $\phi $, $\psi $, $\psi _{k}\in
C_{b}^{1,2}(Q_{T})$,

\begin{description}
\item[(B1)] $G^{\kappa }(p,X,\phi (t,x+\cdot ),\phi (t,x+\cdot ))=G(p,X,\phi
(t,x+\cdot ))$.

\item[(B2)] If $X\geq Y$, $(v-u)(t,\cdot )$ and $(\phi -\psi )(t,\cdot )$
have a global minimum at $x$, then%
\[
G^{\kappa }(p,X,v(t,x+\cdot ),\phi (t,x+\cdot ))\geq G^{\kappa
}(p,Y,u(t,x+\cdot ),\psi (t,x+\cdot )).
\]

\item[(B3)] For $c_{1},c_{2}\in \mathbb{R}$,
\[
G^{\kappa }(p,X,w(t,x+\cdot )+c_{1},\phi (t,x+\cdot
)+c_{2})=G^{\kappa }(p,X,w(t,x+\cdot ),\phi (t,x+\cdot )).
\]

\item[(B4)] If $\psi _{k}(t,\cdot )\rightarrow w(t,\cdot )$ locally uniform
on $\mathbb{R}^{d}$ and $\psi _{k}(t,\cdot )$ uniformly bounded, then $%
G^{\kappa }(p,X,\psi _{k}(t,x+\cdot ),\phi (t,x+\cdot ))\rightarrow
G^{\kappa }(p,X,w(t,x+\cdot ),\phi (t,x+\cdot ))$.
\end{description}

For the equation (\ref{Main-eq}), the corresponding $G^{\kappa }$ is%
\begin{eqnarray*}
G^{\kappa }(p,X,u(t,x+\cdot ),\phi (t,x+\cdot )) =\sup_{(v,q,Q)\in \mathcal{U%
}}\{ \int_{|z|>\kappa }(u(t,x+z)-u(t,x))v(dz) \\
+\int_{|z|\leq \kappa }(\phi (t,x+z)-\phi (t,x))v(dz)+\langle p,q\rangle +%
\frac{1}{2}\text{\textrm{tr}}[XQQ^{T}]\}.
\end{eqnarray*}%
It is easy to check that this $G^{\kappa }$ satisfies the above
assumptions. The above assumptions also hold for the equation
(\ref{Main-geq}).

\begin{remark}
Our assumptions (A1) and (B4) are different from \cite{B-I, J-K}.
This is because that the measures in $\mathcal{U}$ may be singular,
which make the problem more difficulty in taking the limit. See the
following example.
\end{remark}

\begin{example}
Consider $\mathcal{V}=\{ \delta _{x}:x\in (1,2]\}$ and an upper
semicontinuous function $f(x)=I_{[0,1]}(x)$. Let $f_{n}$ be a
sequence of continuous functions such that $f_{n}\rightarrow f$
pointwise. Then it is easy to show that $\sup_{v\in \mathcal{V}}\int
f_{n}(z)v(dz)$ does not tend to $\sup_{v\in \mathcal{V}}\int
f(z)v(dz)$.
\end{example}

\begin{proposition}
Suppose $u\in \mathrm{USC}_{b}(Q_{T})$ ($u\in
\mathrm{LSC}_{b}(Q_{T})$) is a viscosity subsolution (viscosity
supersolution) of (\ref{eq1}). If $u$ is continuous in $x$ and for
$\phi \in C^{1,2}(Q_{T})$, $(t,x)\in Q_{T}$ is a global maximum
point (minimum point) of $u-\phi $, then for each $\kappa \in
(0,1)$ we have%
\[
\partial _{t}\phi (t,x)-G^{\kappa }(D\phi (t,x),D^{2}\phi (t,x),u(t,x+\cdot
),\phi (t,x+\cdot ))\leq 0\ (\geq 0).
\]
\end{proposition}

The proof can be found in \cite{J-K} and the references therein.

In the following, we first extend the matrix lemma in \cite{Cr}.

\begin{theorem}
Suppose that $X,Y\in \mathbb{S}(N)$ satisfy $X\leq Y<\frac{1}{\gamma
}I$ for some $\gamma >0$. Define $X^{\gamma }=X(I-\gamma X)^{-1}$
and $Y^{\gamma
}=Y(I-\gamma Y)^{-1}$. Then $Y^{\gamma }\geq X^{\gamma }\geq X$ and $%
X^{\gamma }\geq -\frac{1}{\gamma }I$.
\end{theorem}

\begin{proof}
For each $A\in \mathbb{S}(N)$ with $A<\frac{1}{\gamma }I$, it is
easy to
check that for each fixed $y\in \mathbb{R}^{N}$,%
\begin{equation}
\max_{x\in \mathbb{R}^{N}}\{ \langle Ax,x\rangle -\frac{1}{\gamma }%
|x-y|^{2}\}=\langle A(I-\gamma A)^{-1}y,y\rangle .  \label{mat}
\end{equation}

From the condition $X\leq Y<\frac{1}{\gamma }I$, we have
\[
\langle Xx,x\rangle -\frac{1}{\gamma }|x-y|^{2}\leq \langle Yx,x\rangle -%
\frac{1}{\gamma }|x-y|^{2}\text{ \ for\ }x,y\in
\mathbb{R}^{N}\text{.}
\]%
Thus for each fixed $y\in \mathbb{R}^{N}$,%
\[
\max_{x\in \mathbb{R}^{N}}\{ \langle Xx,x\rangle -\frac{1}{\gamma }%
|x-y|^{2}\} \leq \max_{x\in \mathbb{R}^{N}}\{ \langle Yx,x\rangle -\frac{1}{%
\gamma }|x-y|^{2}\}.
\]%
By (\ref{mat}) we obtain%
\[
\langle X^{\gamma }y,y\rangle \leq \langle Y^{\gamma }y,y\rangle
\text{ \ for all }y\in \mathbb{R}^{N},
\]%
which yields $Y^{\gamma }\geq X^{\gamma }$. It is easy to check
$X^{\gamma }\geq X$ and $X^{\gamma }\geq -\frac{1}{\gamma }I$. The
proof is complete.
\end{proof}

In particular, we consider%
\[
J_{nd}:=\left(
\begin{array}{cccc}
(n-1)I & -I & \cdots & -I \\
-I & \ddots & \ddots & \vdots \\
\vdots & \ddots & \ddots & -I \\
-I & \cdots & -I & (n-1)I%
\end{array}%
\right) _{nd\times nd}.
\]%
It is easy to prove that $J_{nd}^{2}=nJ_{nd}$ and for each given
$\gamma \in
(0,\frac{1}{n})$, $J_{nd}<\frac{1}{\gamma }I$ and $J_{nd}^{\gamma }=\frac{1}{%
1-n\gamma }J_{nd}$. In the following, we always use $J_{nd}$ for
convenience. Then we immediately get the following corollary.

\begin{corollary}
Let $X_{i}\in \mathbb{S}(d)$, $i=1,\ldots ,n$, satisfy%
\[
\left(
\begin{array}{cccc}
X_{1} & 0 & \cdots & 0 \\
0 & \ddots & \ddots & \vdots \\
\vdots & \ddots & \ddots & 0 \\
0 & \cdots & 0 & X_{n}%
\end{array}%
\right) \leq J_{nd}.
\]%
Then for each given $\gamma \in (0,\frac{1}{n})$, $(I-\gamma
X_{i})^{-1}$
exists for $i=1,\ldots ,n$, and $X_{i}^{\gamma }:=X_{i}(I-\gamma X_{i})^{-1}$%
, $i=1,\ldots ,n$, satisfy $X_{i}^{\gamma }\geq X_{i}$ and%
\[
-\frac{1}{\gamma }I_{nd}\leq \left(
\begin{array}{cccc}
X_{1}^{\gamma } & 0 & \cdots & 0 \\
0 & \ddots & \ddots & \vdots \\
\vdots & \ddots & \ddots & 0 \\
0 & \cdots & 0 & X_{n}^{\gamma }%
\end{array}%
\right) \leq \frac{1}{1-n\gamma }J_{nd}.
\]
\end{corollary}

We now give the main lemma (see Lemma 7.8 in \cite{J-K}).

\begin{lemma}
\label{app-le1}Let $u_{i}\in \mathrm{USC}_{b}(Q_{T})$ be viscosity
subsolutions of%
\[
\partial _{t}u(t,x)-G_{i}(Du(t,x),D^{2}u(t,x),u(t,x+\cdot ))=0,\ i=1,\ldots
,n,
\]%
on $Q_{T}$, where $G_{i}$, $i=1,\ldots ,n$, satisfy (A1)-(A3). Let
$\phi \in C_{b}^{1,2}$ satisfy that $(\bar{t},\bar{x}_{1},\ldots
,\bar{x}_{n})\in
(0,T)\times \mathbb{R}^{nd}$ is a global maximum point of $%
\sum_{i=1}^{n}u_{i}(t,x_{i})-\phi (t,x_{1},\ldots ,x_{n})$.
Moreover, suppose that there exist continuous functions $g_{0}>0$,
$g_{1},\ldots
,g_{n} $ such that%
\[
D^{2}\phi \leq g_{0}(t,x)J_{nd}+\left(
\begin{array}{cccc}
g_{1}(t,x_{1}) & 0 & \cdots & 0 \\
0 & \ddots & \ddots & \vdots \\
\vdots & \ddots & \ddots & 0 \\
0 & \cdots & 0 & g_{n}(t,x_{n})%
\end{array}%
\right) .
\]%
Then for each $\gamma \in (0,\frac{1}{n})$, there exist $b_{i}\in
\mathbb{R}$ and $X_{i}\in \mathbb{S}(d)$, $i=1,\ldots ,n$, such that

\begin{description}
\item[(i)] $b_{1}+\cdots +b_{n}=\partial _{t}\phi (\bar{t},\bar{x}%
_{1},\ldots ,\bar{x}_{n})$;

\item[(ii)]
\[
-\frac{g_{0}(\bar{t},\bar{x})}{\gamma }I_{nd}\leq \left(
\begin{array}{cccc}
X_{1} & 0 & \cdots & 0 \\
0 & \ddots & \ddots & \vdots \\
\vdots & \ddots & \ddots & 0 \\
0 & \cdots & 0 & X_{n}%
\end{array}%
\right) -K_{nd}\leq \frac{g_{0}(\bar{t},\bar{x})}{1-n\gamma }J_{nd},
\]%
where $\bar{x}=(\bar{x}_{1},\ldots ,\bar{x}_{n})$ and%
\[
K_{nd}=\left(
\begin{array}{cccc}
g_{1}(\bar{t},\bar{x}_{1}) & 0 & \cdots & 0 \\
0 & \ddots & \ddots & \vdots \\
\vdots & \ddots & \ddots & 0 \\
0 & \cdots & 0 & g_{n}(\bar{t},\bar{x}_{n})%
\end{array}%
\right) ;
\]

\item[(iii)] for each $i=1,\ldots ,n$,%
\[
b_{i}-G_{i}(D_{x_{i}}\phi (\bar{t},\bar{x}),X_{i},\phi (\bar{t},\bar{x}%
_{1},\ldots ,\bar{x}_{i-1},\bar{x}_{i}+\cdot ,\bar{x}_{i+1},\ldots ,\bar{x}%
_{n}))\leq 0.
\]
\end{description}
\end{lemma}

\begin{remark}
Applying the above matrix inequalities, the proof in \cite{J-K}
still holds.
\end{remark}

\begin{remark}
If $u_{i}$ is continuous in $x$, we can further get that for each
$\kappa
\in (0,1)$,%
\[
b_{i}-G_{i}^{\kappa }(D_{x_{i}}\phi (\bar{t},\bar{x}),X_{i},u_{i}(\bar{t},%
\bar{x}_{i}+\cdot ),\phi (\bar{t},\bar{x}_{1},\ldots ,\bar{x}_{i-1},\bar{x}%
_{i}+\cdot ,\bar{x}_{i+1},\ldots ,\bar{x}_{n}))\leq 0.
\]
\end{remark}

We now give the main theorem, which combines the methods in
\cite{Peng2008} and \cite{B-I}.

\begin{theorem}
\label{app-th1}(Domination theorem) Let $u_{i}\in \mathrm{USC}%
_{b}([0,T)\times \mathbb{R}^{d})$ be viscosity subsolutions of%
\[
\partial _{t}u(t,x)-G_{i}(Du(t,x),D^{2}u(t,x),u(t,x+\cdot ))=0,\ i=1,\ldots
,n,
\]%
on $Q_{T}$, where $G_{i}$ and $G_{i}^{\kappa }$, $i=1,\ldots ,n$,
satisfy (A1)-(A3) and (B1)-(B4). We suppose also that

\begin{description}
\item[(i)] There exists a constant $C>0$ such that for each $\kappa \in
(0,1) $, $p,q\in \mathbb{R}^{d}$, $X$, $Y\in \mathbb{S}(d)$, $u\in \mathrm{SC%
}_{b}(Q_{T})$, $\varphi \in C^{1,2}(Q_{T})$, $\psi _{1}\in C_{b}^{1}(\mathbb{%
R}^{d})$ and $\psi _{2}\in C_{b}^{2}(\mathbb{R}^{d})$,
\begin{eqnarray*}
&&|G_{1}^{\kappa }(p,X,u(t,\cdot )+\psi _{1}(\cdot ),\varphi
(t,\cdot )+\psi
_{2}(\cdot ))-G_{1}^{\kappa }(q,Y,u(t,\cdot ),\varphi (t,\cdot ))| \\
&\leq &C(|p-q|+||X-Y||+\sup_{x\in \mathbb{R}^{d}}|D\psi
_{1}(x)|+\sup_{x\in \mathbb{R}^{d}}(|D\psi _{2}(x)|+|D^{2}\psi
_{2}(x)|));
\end{eqnarray*}

\item[(ii)] For given constants $\beta _{i}>0$, $i=1,\ldots ,n$, the
following domination condition holds for $G_{i}$: for each
$(t,x_{1},\ldots
,x_{n})\in (0,T)\times \mathbb{R}^{nd}$, $(p_{i},X_{i})\in \mathbb{R}%
^{d}\times \mathbb{S}(d)$ and $\phi _{i}\in C^{1,2}(Q_{T})$ such that $%
\sum_{i=1}^{n}\beta _{i}p_{i}=0$, $\sum_{i=1}^{n}\beta _{i}X_{i}\leq 0$, $%
\sum_{i=1}^{n}\beta _{i}u_{i}(t,x_{i}+\cdot )\leq
\sum_{i=1}^{n}\beta
_{i}u_{i}(t,x_{i})$ and $\sum_{i=1}^{n}\beta _{i}D\phi _{i}(t,x_{i})=0$,%
\[
\sum_{i=1}^{n}\beta _{i}G_{i}^{\kappa
}(p_{i},X_{i},u_{i}(t,x_{i}+\cdot ),\phi _{i}(t,x_{i}+\cdot ))\leq
o_{\kappa }(1)\text{ as }\kappa \rightarrow 0.
\]
\end{description}

Then a similar domination also holds for the solutions: if
$u_{i}(0,\cdot )\in $ $C_{b.Lip}(\mathbb{R}^{d})$ and
$\sum_{i=1}^{n}\beta _{i}u_{i}(0,\cdot )\leq 0$, then
$\sum_{i=1}^{n}\beta _{i}u_{i}(t,\cdot )\leq 0$ for all $t>0$.
\end{theorem}

\begin{proof}
For each given $\bar{\delta}>0$, it is easy to check that for each
$1\leq i\leq n$, $\tilde{u}_{i}:=u_{i}-\bar{\delta}/(T-t)$ is a
viscosity solution
of%
\begin{equation}
\partial _{t}\tilde{u}(t,x)-G_{i}(D\tilde{u}(t,x),D^{2}\tilde{u}(t,x),\tilde{%
u}(t,x+\cdot ))\leq -c,\ c:=\bar{\delta}/T^{2}.  \label{eq2}
\end{equation}%
For each $\lambda >0$, the sup-convolution of $\tilde{u}_{i}$ is defined by%
\[
\tilde{u}_{i}^{\lambda }(t,x):=\sup_{y\in \mathbb{R}^{d}}\{ \tilde{u}%
_{i}(t,y)-\frac{|x-y|^{2}}{\lambda }\}=\sup_{y\in \mathbb{R}%
^{d}}\{u_{i}(t,y)-\frac{|x-y|^{2}}{\lambda
}\}-\frac{\bar{\delta}}{T-t}.
\]%
The function $\tilde{u}_{i}^{\lambda }$ is upper semicontinuous in
$(t,x)$ and continuous in $x$. Moreover, $\tilde{u}_{i}^{\lambda }$
is still a
viscosity solution of (\ref{eq2}) (see Lemma 7.3 in \cite{J-K}). Note that $%
u_{i}(0,\cdot )\in $ $C_{b.Lip}(\mathbb{R}^{d})$, then we can choose
a small
$\lambda _{0}>0$ such that for each $\lambda \leq \lambda _{0}$, $%
\sum_{i=1}^{n}\beta _{i}\tilde{u}_{i}^{\lambda }(0,\cdot )\leq 0$. Since $%
\sum_{i=1}^{n}\beta _{i}u_{i}\leq 0$ follows from $\sum_{i=1}^{n}\beta _{i}%
\tilde{u}_{i}\leq 0$ in the limit $\bar{\delta}\downarrow 0$ and $%
\sum_{i=1}^{n}\beta _{i}\tilde{u}_{i}^{\lambda }\downarrow
\sum_{i=1}^{n}\beta _{i}\tilde{u}_{i}$ as $\lambda \downarrow 0$, it
suffices to prove the theorem under the additional assumptions:
$\hat{u}_{i}$ is a viscosity solution of (\ref{eq2}), $\hat{u}_{i}$
is continuous in $x$
and $\lim_{t\rightarrow T}\hat{u}_{i}(t,x)=-\infty $ uniformly in $%
[0,T)\times \mathbb{R}^{d}$. To prove the theorem, we assume to the
contrary
that%
\[
\sup_{(t,x)\in \lbrack 0,T)\times \mathbb{R}^{d}}\sum_{i=1}^{n}\beta _{i}%
\hat{u}_{i}(t,x)=m_{0}>0.
\]%
We define for $x=(x_{1},x_{2},\ldots ,x_{n})\in \mathbb{R}^{nd}$%
\[
\phi _{\varepsilon ,\beta }(x)=\frac{1}{2\varepsilon }\sum_{1\leq
i<j\leq n}|x_{i}-x_{j}|^{2}+\psi _{\beta }(x_{1}),
\]%
where $\psi _{\beta }(x_{1}):=\psi (\beta x_{1})$ and $\psi \in C_{b}^{2}(%
\mathbb{R}^{d})$ such that $\psi (x_{1})=0$ for $|x_{1}|\leq 1$ and
$\psi (x_{1})=2\sum_{i=1}^{n}\beta _{i}||u_{i}||_{\infty }$ for
$|x_{1}|\geq 2$. It is easy to check that $\psi _{\beta
}(x_{1})=2\sum_{i=1}^{n}\beta _{i}||u_{i}||_{\infty }$ for
$|x_{1}|\geq 2/\beta $ and $\sup \{|D\psi
_{\beta }|+||D^{2}\psi _{\beta }||\} \rightarrow 0$ as $\beta \rightarrow 0$%
. We define%
\[
M_{\varepsilon ,\beta }=\sup_{(t,x)\in \lbrack 0,T)\times
\mathbb{R}^{nd}}\{ \sum_{i=1}^{n}\beta _{i}\hat{u}_{i}(t,x_{i})-\phi
_{\varepsilon ,\beta }(x)\}.
\]%
By the construction of $\psi _{\beta }$, for small $\varepsilon $, $\beta >0$%
, the maximum of the above function is achieved at some $(\bar{t},\bar{x})=(%
\bar{t},\bar{x}_{1},\ldots ,\bar{x}_{n})$ with $|\bar{x}_{1}|\leq
2/\beta $, $\sum_{1\leq i<j\leq n}|\bar{x}_{i}-\bar{x}_{j}|^{2}\leq
(2\varepsilon \sum_{i=1}^{n}\beta _{i}||u_{i}||_{\infty })^{1/2}$
and $M_{\varepsilon ,\beta }\geq m_{0}/2>0$. Thus there exists a
constant $T_{0}<T$ independent of $\varepsilon $ and $\beta $ such
that $\bar{t}\leq T_{0}$. For fixed small $\beta $, we can check
that (see Lemma 3.1 in \cite{CIL})

\begin{description}
\item[1)] $\frac{1}{\varepsilon }\sum_{1\leq i<j\leq n}|\bar{x}_{i}-\bar{x}%
_{j}|^{2}\rightarrow 0$ as $\varepsilon \rightarrow 0$.

\item[2)] $\lim_{\varepsilon \rightarrow 0}M_{\varepsilon ,\beta
}=\sup_{(t,z)\in \lbrack 0,T)\times \mathbb{R}^{d}}\{
\sum_{i=1}^{n}\beta
_{i}\hat{u}_{i}(t,z)-\psi _{\beta }(z)\}=\sum_{i=1}^{n}\beta _{i}\hat{u}_{i}(%
\hat{t},\hat{z})-\psi _{\beta }(\hat{z})\geq m_{0}/2$, where $(\hat{t},\hat{z%
})$ is any limit point of $(\bar{t},\bar{x}_{1})$.
\end{description}

Since $\hat{u}_{i}\in \mathrm{USC}([0,T)\times \mathbb{R}^{d})$ and $%
\sum_{i=1}^{n}\beta _{i}\hat{u}_{i}(0,\cdot )\leq 0$, it is easy to get $%
\hat{t}>0$. Thus $\bar{t}$ must be strictly positive for small $\varepsilon $%
. Applying Lemma \ref{app-le1} at the point $(\bar{t},\bar{x})=(\bar{t},\bar{%
x}_{1},\ldots ,\bar{x}_{n})$ and taking $\gamma =1/(2n)$, we obtain
that
there exist $b_{i}\in \mathbb{R}$ and $X_{i}\in \mathbb{S}(d)$ for $%
i=1,\ldots ,n$ such that $\sum_{i=1}^{n}\beta _{i}b_{i}=0$,
\[
\left(
\begin{array}{cccc}
\beta _{1}X_{1} & 0 & \cdots & 0 \\
0 & \ddots & \ddots & \vdots \\
\vdots & \ddots & \ddots & 0 \\
0 & \cdots & 0 & \beta _{n}X_{n}%
\end{array}%
\right) -\left(
\begin{array}{cccc}
D^{2}\psi _{\beta }(\bar{x}_{1}) & 0 & \cdots & 0 \\
0 & 0 & \ddots & \vdots \\
\vdots & \ddots & \ddots & 0 \\
0 & \cdots & 0 & 0%
\end{array}%
\right) \leq \frac{2}{\varepsilon }J_{nd},
\]%
and for each $\kappa \in (0,1)$,%
\[
b_{i}-G_{i}^{\kappa
}(p_{i},X_{i},\hat{u}_{i}(\bar{t},\bar{x}_{i}+\cdot
),\beta _{i}^{-1}\phi _{\varepsilon ,\beta }(\bar{x}_{1},\ldots ,\bar{x}%
_{i-1},\bar{x}_{i}+\cdot ,\bar{x}_{i+1},\ldots ,\bar{x}_{n}))\leq
-c,
\]%
where $p_{i}=\beta _{i}^{-1}D_{x_{i}}\phi _{\varepsilon ,\beta
}(\bar{x})$.
It is easy to check that $\sum_{i=1}^{n}\beta _{i}p_{i}=D\psi _{\beta }(\bar{%
x}_{1})$ and $\sum_{i=1}^{n}\beta _{i}X_{i}\leq D^{2}\psi _{\beta }(\bar{x}%
_{1})$. Thus for each $\kappa \in (0,1)$,
\begin{eqnarray*}
&&-c\sum_{i=1}^{n}\beta _{i}=-\sum_{i=1}^{n}\beta
_{i}b_{i}-c\sum_{i=1}^{n}\beta _{i} \\
&\geq &-\sum_{i=1}^{n}\beta _{i}G_{i}^{\kappa }(p_{i},X_{i},\hat{u}_{i}(\bar{%
t},\bar{x}_{i}+\cdot ),\beta _{i}^{-1}\phi _{\varepsilon ,\beta }(\bar{x}%
_{1},\ldots ,\bar{x}_{i-1},\bar{x}_{i}+\cdot ,\bar{x}_{i+1},\ldots ,\bar{x}%
_{n})) \\
&\geq &o_{\kappa }(1)-\beta _{1}G_{1}^{\kappa }(p_{1},X_{1},\hat{u}_{1}(\bar{%
t},\bar{x}_{1}+\cdot ),\beta _{1}^{-1}\phi _{\varepsilon ,\beta }(\bar{x}%
_{1}+\cdot ,\tilde{x})) \\
&&+\beta _{1}G_{1}^{\kappa }(p_{1}-\tilde{p}_{1},X_{1}-\tilde{X}_{1},\hat{u}%
_{1}(\bar{t},\bar{x}_{1}+\cdot )-l(\cdot ),\beta _{1}^{-1}\phi
_{\varepsilon
,\beta }(\bar{x}_{1}+\cdot ,\tilde{x})-l(\cdot )) \\
&\geq &o_{\kappa }(1)-C\sup \{|D\psi _{\beta }|+||D^{2}\psi _{\beta
}||\},
\end{eqnarray*}%
where $\tilde{x}=(\bar{x}_{2},\ldots ,\bar{x}_{n})$,
$\tilde{p}_{1}=\beta _{1}^{-1}D\psi _{\beta }(\bar{x}_{1})$,
$\tilde{X}_{1}=\beta _{1}^{-1}D^{2}\psi _{\beta }(\bar{x}_{1})$ and
$l(\cdot )=\beta _{1}^{-1}\psi _{\beta }(\bar{x}_{1}+\cdot )$. The
right side tends to zero as $\kappa \rightarrow 0$ first and then
$\beta \rightarrow 0$, which induces a contradiction. The proof is
complete.
\end{proof}

\begin{remark}
The above theorem still holds for general $G$, which may contain $%
(t,x,u,Du,D^{2}u,u(t,\cdot ))$.
\end{remark}

We have the following corollaries which are important in this paper.

\begin{corollary}
\label{app-co1}(Comparison theorem) Let the functions $u\in \mathrm{USC}%
_{b}([0,T)\times \mathbb{R}^{d})$ and $v\in
\mathrm{LSC}_{b}([0,T)\times \mathbb{R}^{d})$ be respectively a
viscosity subsolution and a viscosity supersolution of (\ref{eq1}).
Suppose $G$ and $G^{\kappa }$ satisfy (A1)-(A3) and (B1)-(B4) and
the condition (i) of Theorem \ref{app-th1}. Furthermore, we suppose
that

\begin{description}
\item[(iii)] For each $(p,X)\in \mathbb{R}^{d}\times \mathbb{S}(d)$, $u\in
\mathrm{SC}_{b}(Q_{T})$, $\phi _{1}$, $\phi _{2}\in C^{1,2}(Q_{T})$ with $%
D\phi _{1}(t,x)=D\phi _{2}(t,x)$, as $\kappa \rightarrow 0$,
\[
|G^{\kappa }(p,X,u(t,x+\cdot ),\phi _{1}(t,x+\cdot ))-G^{\kappa
}(p,X,u(t,x+\cdot ),\phi _{2}(t,x+\cdot ))|\rightarrow 0.
\]
\end{description}

If $u(0,\cdot )$, $v(0,\cdot )\in C_{b.Lip}(\mathbb{R}^{d})$ and
$u(0,\cdot )\leq v(0,\cdot )$, then $u(t,\cdot )\leq v(t,\cdot )$
for all $t>0$.
\end{corollary}

\begin{corollary}
(Sub-additivity) Let $u^{\varphi }$ denote the viscosity solution of (\ref%
{eq1}) with initial condition $\varphi $. Suppose $G$ and $G^{\kappa
}$ satisfy all the conditions of Corollary \ref{app-co1} and the
following condition:

\begin{description}
\item[(iv)] For each $(p_{i},X_{i})\in \mathbb{R}^{d}\times \mathbb{S}(d)$, $%
u_{i}\in \mathrm{SC}_{b}(Q_{T})$, $\phi _{i}\in C^{1,2}(Q_{T})$, $i=1,2$,%
\begin{eqnarray*}
&&G^{\kappa }(p_{1}+p_{2},X_{1}+X_{2},(u_{1}+u_{2})(t,\cdot ),(\phi
_{1}+\phi _{2})(t,\cdot )) \\
&\leq &G^{\kappa }(p_{1},X_{1},u_{1}(t,\cdot ),\phi _{1}(t,\cdot
))+G^{\kappa }(p_{2},X_{2},u_{2}(t,\cdot ),\phi _{2}(t,\cdot )).
\end{eqnarray*}
\end{description}

Then $u^{\varphi +\psi }\leq u^{\varphi }+u^{\psi }$ for each $\varphi $, $%
\psi \in C_{b.Lip}(\mathbb{R}^{d})$.
\end{corollary}

\begin{corollary}
(Convexity) Let $u^{\varphi }$ denote the viscosity solution of
(\ref{eq1}) with initial condition $\varphi $. Suppose $G$ and
$G^{\kappa }$ satisfy all the conditions of Corollary \ref{app-co1}
and the following condition:

\begin{description}
\item[(v)] For each $\lambda \in (0,1)$, $(p_{i},X_{i})\in \mathbb{R}%
^{d}\times \mathbb{S}(d)$, $u_{i}\in \mathrm{SC}_{b}(Q_{T})$, $\phi
_{i}\in
C^{1,2}(Q_{T})$, $i=1,2$,%
\begin{eqnarray*}
&&G^{\kappa }(\lambda (p_{1},X_{1},u_{1}(t,\cdot ),\phi _{1}(t,\cdot
))+(1-\lambda )(p_{2},X_{2},u_{2}(t,\cdot ),\phi _{2}(t,\cdot ))) \\
&\leq &\lambda G^{\kappa }(p_{1},X_{1},u_{1}(t,\cdot ),\phi
_{1}(t,\cdot ))+(1-\lambda )G^{\kappa }(p_{2},X_{2},u_{2}(t,\cdot
),\phi _{2}(t,\cdot )).
\end{eqnarray*}
\end{description}

Then $u^{\lambda \varphi +(1-\lambda )\psi }\leq \lambda u^{\varphi
}+(1-\lambda )u^{\psi }$ for each $\varphi $, $\psi \in C_{b.Lip}(\mathbb{R}%
^{d})$.
\end{corollary}

For our main equations (\ref{Main-eq}) and (\ref{Main-geq}), it is
easy to check that all assumptions (i)-(v) hold. Perron's method for
(\ref{eq1}) still holds (the proof is similar to Proposition 1 in
\cite{A-T}).

For each $\varphi \in C_{b}^{2}(\mathbb{R}^{d})$, it is easy to find
constants $M_{1}$ and $M_{2}$ such that $u(t,x):=M_{1}t+\varphi $ and $%
v(t,x):=M_{2}t+\varphi $ are respectively the viscosity subsolution
and supersolution of (\ref{Main-eq}) and (\ref{Main-geq}). By
Perron's method and approximation, the solutions of (\ref{Main-eq})
and (\ref{Main-geq}) exist for each $\varphi \in
C_{b.Lip}(\mathbb{R}^{d})$.

\end{document}